\documentclass[preprint,3p]{elsarticle}

\usepackage{amsmath,amssymb,latexsym,color}
\usepackage{amssymb}
 \usepackage{amsthm}

\newtheorem{coro}{Corollary}
\newtheorem{example}{Example}
\newtheorem{definition}{Definition}
\newtheorem{theorem}{Theorem}
\newtheorem{lemma}[theorem]{Lemma}
\newtheorem{note}{Remark}

\newcommand{\CC}{\mathbb{C}}

\newcommand{\NN}{\mathbb{N}}
\newcommand{\PP}{\mathbb{P}}
\newcommand{\RR}{\mathbb{R}}
\newcommand{\ZZ}{\mathbb{Z}}
\newcommand{\LL}{\mathcal{L}}

\newcommand{\dsty}{\displaystyle}

\begin{document}

\begin{frontmatter}


%% use the tnoteref command within \title for footnotes; %% use the tnotetext command for the associated footnote; %% use the fnref command within \author or \address for footnotes; %% use the fntext command for the associated footnote; %% use the corref command within \author for corresponding author footnotes; %% use the cortext command for the associated footnote; %% use the ead command for the email address, %% and the form \ead[url] for the home page: %%

% \title{}
%% \tnotetext[label1]{} %% \author{Name\corref{cor1}\fnref{label2}} %% \ead{email address} %% \ead[url]{home page} %% \fntext[label2]{} %% \cortext[cor1]{} %% \address{Address\fnref{label3}} %% \fntext[label3]{}

\title{On orthogonal polynomials with respect to  a class of differential operators.}

\author{Jorge A. Borrego–Morell}
%\address[label1]{Universidad Carlos III de Madrid\\
%Departamento de Matem\'aticas, Avenida de la Universidad, 30 CP-28911, Legan\'es, Madrid.} \ead{jbmorell@gmail.com}
%
%\fntext[grant]{Research partially supported by Direcci\'on General de Investigaci\'on, Ministerio de Ciencia e
%Innovaci\'on of Spain, under grant MTM2009-12740-C03-01}

\begin{abstract}
We consider orthogonal polynomials with respect to a linear differential operator

$$
\LL^{(M)}=\sum_{k=0}^{M}\rho_{k}(z)\frac{d^k}{dz^k},
$$

\noindent where  $\{\rho_k\}_{k=0}^{M}$ are complex polynomials such that $deg[\rho_k]\leq k, 0\leq k \leq
M$, with equality for at least one index. We  analyze the uniqueness and zero location of these polynomials. An
interesting phenomenon  occurring in this kind of orthogonality is the existence of operators for which the associated
sequence of orthogonal polynomials reduces to a finite set. For a given operator, we find a classification  of the
measures for which it is possible to guarantee the existence of an infinite sequence of orthogonal polynomials, in
terms of a linear system of difference equations with varying coefficients.  Also, for the case of a first-order
differential operator, we  locate  the zeros and establish the strong asymptotic behavior of these polynomials.

\end{abstract}

\begin{keyword}
Orthogonal polynomials \sep linear differential operators \sep zero location \sep asymptotic behavior.

\MSC[2010] 42C05 \sep 47E05.

%% MSC codes here, in the form: \MSC code \sep code %% or \MSC[2008] code \sep code (2000 is the default)

\end{keyword}

\end{frontmatter}

\section{Introduction}

Orthogonality with respect to a linear homogeneous differential operator was introduced in \cite{ApLoMa02}  as a
 generalization of the notion of orthogonal polynomials. There, the authors show that the concept  of Chebyshev system
plays a fundamental role to solve the problem of the uniqueness of the sequence of the polynomials. A further study of some algebraic and analytic properties of this type of orthogonality is done in \cite{BePiMaUr11,BP12} for
some first and second--order linear homogeneous differential operators. Formally, orthogonality  with respect to a linear homogeneous differential operator is defined as follows,
\begin{definition}\label{DoD}
Assume that $\mu$ is a finite positive Borel measure on the real line  and let $\{\rho_k\}_{k=0}^{M}$ be a
set of functions such that,
$$\int |x^j\rho_k(x)|d\mu(x)<\infty, \quad 0\leq j <\infty,$$
\noindent for all $k=0,\ldots,M$.  Denote
\begin{eqnarray}\label{gop}
\LL^{(M)}=\sum_{k=0}^{M}\rho_{k}(x)\frac{d^{k}}{dx^{k}},
\end{eqnarray}
an operator acting over the space of polynomials $\PP$.

We say that $\{Q_{n}\}_{n=0}^{\infty}$,  is a sequence of orthogonal polynomials with respect to the pair
$(\LL^{(M)},\mu)$ if $deg[Q_n]\leq n$ and

\begin{eqnarray}\label{GO}
\dsty \int \LL^{(M)}[Q_n(x)]P(x)d\mu(x)=0,
\end{eqnarray}

\noindent for any polynomial $P$ such that $deg[P]\leq n-1$.

 \end{definition}

We recall that in the definition of orthogonality with respect to a differential operator given in \cite{ApLoMa02},
$\rho_M$ is assumed to be equal to $1$, but we shall drop this assumption. The determination of the sequence of these
polynomials can be reduced to the solution of a system of $n$ algebraic linear homogeneous  equations on the $n+1$
coefficients of $Q_n$, thus the existence is guaranteed. Unlike systems of orthogonal polynomials, it is not possible
to affirm uniqueness up to a constant factor and this turns out to be  a difficult problem. We say that an
index $n$ is normal if for this $n$ the solution is uniquely determined up to a constant multiplicative factor.  For a fixed non--negative integer $n$, $Q_n$ will be referred to as the orthogonal polynomial with respect to the pair
$(\LL^{(M)},\mu)$ associated with the index $n$, which in general is not necessarily  unique.

In this  manuscript, we show  necessary and sufficient conditions for the normality of an index $n\in\ZZ_+$ for exactly
solvable differential operators and study   properties of the sequence of orthogonal polynomials. Exactly
solvable operators arise from quantum mechanics and were introduced by Turbiner in \cite{turb92a,turb92b},
\begin{definition}\label{exactlysolv} 
 Let  $\mathbb{P}_n$ be  the space of all polynomials of degree at most $n$. A linear differential operator of  $M$-th order $\LL^{(M)}$ is called 
exactly-solvable if   $\LL^{(M)}[\mathbb{P}_n]\subseteq \mathbb{P}_n$, for all $n\geq 0$, with equality for at least one index $n$.
\end{definition}
Notice that any exactly-solvable operator has the form
\begin{eqnarray*}
\LL^{(M)}=\sum_{k=0}^{M}\rho_{k}(x)\frac{d^{k}}{dx^{k}},
\end{eqnarray*}
 where   $\dsty \rho_{k}(x)=\sum_{j=0}^{k}\rho_{k,j}x^j$ are polynomials  which satisfy the condition $deg[\rho_{k}]\leq k$ with equality for at least one index $k$.  

Some of the techniques used here could be extended to some degree to the  general case of linear homogeneous
differential  operators with polynomial coefficients including the  class of Heine--Stieltjes operators
as well as the lowering and raising  operators with polynomial coefficients, but we will not dwell on it.

The paper is organized as follows. In Section \ref{examples} we present connections between  this type of orthogonality and
some inner products and classify  the exactly solvable  operators for which this concept of
orthogonality reduces to an inner product. In Section \ref{exists} we  give necessary and sufficient conditions
for the normality of an index $n$. The analysis  of the existence of infinite sequences of polynomials
$\{Q_n\}_{n=m+1}^{\infty}$, for some positive $m$, with $deg[Q_n]=n$, in terms of a
linear system of difference equations with varying coefficients is done in  section \ref{Infinite}. In Section \ref{Zeroloc},
we study the location of the zeros for the  polynomials $Q_n$, and  in Section \ref{Polar}, for a first-order
differential operator, we obtain a curve that contains  the accumulation points of the zeros of the polynomials  giving also the
strong asymptotic behavior of the polynomials. Finally, in Section \ref{Cr}, we discuss some possible extensions of the results.

\section{Applications and relation to some  inner products}\label{examples}

Let us see some examples where orthogonality with respect to a differential operator reduces in some sense to orthogonality with respect to an inner product.

\begin{enumerate}

\item When $M=0$, we obtain the classical construction of orthogonal polynomials with respect to a standard inner
    product

$$\int Q_n(x)P(x)d\mu(x)=0, \quad deg[P]\leq n-1.$$

\item Let $\zeta \in \CC$ be fixed and consider the differential operator $\LL_{\zeta}:W^{1,2}(\mu)\rightarrow
    L^{2}(\mu)$

$$\LL_{\zeta}[f](x)=f(x)+(x-\zeta)f^{\prime}(x),$$
where $W^{1,2}(\mu)=\{f\in L^{2}(\mu): f^{\prime}\in L^{2}(\mu)\}$  is the Sobolev space of index 1.
 Let us consider a  positive Borel measure $\mu$
 supported on a subset $\Delta \subset \RR$. The \emph{polar polynomial} associated to $\mu$, see \cite{PiBeUr10},
 is defined as the polynomial $Q_n$ of degree $n$ orthogonal with respect to $(\LL_{\zeta},\mu)$. Let us consider

\begin{eqnarray*}
\Pi_{0,\zeta}&=&1, \\ \Pi_{n+1,\zeta}(z)&=&(z-\zeta) Q_n(z), \quad n\geq 0.
\end{eqnarray*}

Then it is not difficult to see that the family of polynomials $\{ \Pi_{n,\zeta}\}_{n=0}^{\infty}$  is orthogonal
with  respect to the Sobolev inner product
\begin{equation}\nonumber
\langle f,g\rangle_{\zeta}= \eta f(\zeta) g(\zeta)+ \int_{\Delta} f^{\prime}(x)\, g^{\prime}(x)\, d\mu(x),
\end{equation}
for some  $\eta>0$.
\end{enumerate}

The authors in \cite{BePiMaUr11} give a detailed study of this family for the case in which $\mu=\mu_{\lambda}$,
$\lambda >-\frac{1}{2}$, is the (classical) Gegenbauer or ultraspherical measure, i.e.
$d\mu_{\lambda}(x)=(1-x^2)^{\lambda -\frac{1}{2}} dx$.

Before we state our first result we remind that a Bochner--Krall type operator is an  operator for
which the associated sequence of eigenpolynomials is orthogonal with respect to some finite positive Borel measure, see \cite{bochner29,krall38,krall40}. Given a
finite positive Borel measure on the real line $\mu$, consider the moment functional $\sigma:\PP\longrightarrow
\RR$ defined by

\begin{eqnarray}\label{functdef}
\dsty \langle\sigma,x^n\rangle=\int x^n d\mu(x).
\end{eqnarray}

Following \cite{kly94}, we define the moment functional  $\sigma^{\prime}$, the derivative of $\sigma$, and
$\phi\sigma$, the multiplication of $\sigma$ times a polynomial $\phi$ as:

\begin{eqnarray}\label{rules}
\langle\sigma^{\prime},\psi\rangle&=&-\langle\sigma,\psi^{\prime}\rangle,\\ \nonumber
\langle\phi\sigma,\psi\rangle&=&\langle\sigma,\phi\psi\rangle,
\end{eqnarray}

\noindent for all $\psi\in \PP$.

Consider now  the bilinear form on $\PP$

\begin{equation}\label{relortho}
\dsty [Q,P] =\int \LL^{(M)}[Q(x)]P(x)d\mu(x).
\end{equation}

In general, it is not possible to affirm that  the bilinear form \eqref{relortho}  defines an inner product.  The
following theorem characterizes the exactly solvable operators and the measures $\mu$ for which  \eqref{relortho}  defines an inner product.

\begin{theorem}\label{innerprod}
Let  $\LL^{(M)}$ be an exactly  solvable and $\mu$ a positive Borel measure. Then, a necessary and sufficient
condition for \eqref{relortho} to be an inner product is that:
\begin{itemize}
\item  $\LL^{(M)}$ is a Bochner-Krall operator and $\mu$ is the measure such that the polynomial eigenfunctions of the operator form a system of orthogonal polynomials.
\item  $\LL^{(M)}$ has positive eigenvalues.
\end{itemize}

In such a case, the sequence of monic polynomials $\{Q_n\}_{n=0}^{\infty}$ orthogonal with respect to $(\LL^{(M)},\mu)$ coincide
with the monic orthogonal polynomials $\{P_n\}_{n=0}^{\infty}$ with respect to the measure $\mu$.

\end{theorem}

\proof

Suppose that the relation \eqref{relortho} defines  an inner product. We have then that \eqref{relortho} is
symmetric, taking into account \eqref{functdef} and \eqref{rules} we have,

\begin{eqnarray*}
\dsty \int \LL^{(M)}[Q](x)P(x)d\mu(x)=\langle\LL^{(M)}[Q](x)\sigma,P(x)\rangle=[ Q,P]\\=
[P,Q]=\int \LL^{(M)}[P](x)Q(x)d\mu(x)=\langle\LL^{(M)}[P]\sigma,Q\rangle;
\end{eqnarray*}

\noindent that is,

$$\langle\LL^{(M)}[Q]\sigma,P\rangle=\langle\LL^{(M)}[P]\sigma,Q\rangle,$$

\noindent and from  \cite[Theorem 2.4( ii) implies i)]{kly94} we have that $\LL^{(M)}$ is a Bochner-Krall operator and $\mu$ is the measure with respect  to which the polynomial eigenfunctions of the operator form a system of orthogonal polynomials.

The second condition follows from the fact that $[.,.]$ must define a positive definite bilinear form.

The converse implication is straightforward.

The assertion  that the sequence of monic polynomials $\{Q_n\}_{n=0}^{\infty}$ coincides with the monic orthogonal polynomials
$\{P_n\}_{n=0}^{\infty}$ with respect to the measure $\mu$ follows from the fact that $\LL^{(M)}$ is an exactly solvable
operator with positive eigenvalues $\{\lambda_n\}_{n=0}^{\infty}$ which implies that condition \eqref{GO} is equivalent
to solving

$$\LL^{(M)}[Q_n]=\lambda_nP_n,$$

\noindent from this  we deduce that $Q_n=P_n$. \qed

We also mention that for the case of first  and second order differential operators, the $n$--th orthogonal
polynomial associated to an index $n$   can be interpreted as the the equilibrium points of a flux of  a complex
potential due to a system of fixed points, cf. \cite{PiBeUr10,BP12}.

\section{Necessary and sufficient conditions for the normality of an index}\label{exists}

In this  section we give necessary and sufficient conditions for the normality of an index $n$ for the class  of
exactly solvable operators. As $\dsty \LL^{(M)}=\sum_{k=0}^{M}\rho_{k}(x)\frac{d^{k}}{dx^{k}}$ is an exactly solvable
operator, following  \eqref{GO},  it is not difficult to see that the monic orthogonal polynomials with respect to the pair
$(\LL^{(M)},\mu)$ associated to an index $n$ are linear combinations of a monic polynomial solution of
\begin{eqnarray}
\label{SolSet1} \LL^{(M)}[y]&=&\lambda_nP_n,
\end{eqnarray}
\noindent and a monic polynomial solution of
\begin{eqnarray}
\label{SolSet2} \LL^{(M)}[y]&=&0,
\end{eqnarray}
here  $ P_n$ denotes the $n$-th monic orthogonal polynomials with respect to $\mu$ and $\dsty
\lambda_n=\sum_{k=0}^{M}\rho_{k,k}\frac{n!}{(n-k)!}$ is  the coefficient associated to the factor $x^n$ of the polynomial $\LL^{(M)}[x^n]$.  By convention, $\dsty \frac{n!}{(n-k)!}=0$ when $k>n$. In the sequel we shall assume that $\lambda_n$  will denote this coefficient.

Before we state the results of the section  we show with an example that in general we do not have normality of an index for the class of
operators that we consider.

\begin{example}\label{counterex}[Second order differential operator]
Suppose that $M=2$ and define $\LL[f]=f^{\prime \prime} -2xf^{\prime}+2f, f \in \PP$. Notice that the eigenfunctions of this operator are the Hermite polynomials
$\{H_n\}_{n=0}^{\infty}$ with eigenvalues $\lambda_n=2(1-n)$ and that $\LL[x]=0$. Consider the measure  $\dsty
d\mu(x)=\frac{e^{-x^2}dx}{x^2+1}$ supported on $\RR$ and denote by   $\{P_n\}_{n=0}^{\infty}$ the sequence of monic
orthogonal polynomials with respect to $\mu$. Notice that  if $n>3$ the polynomial $P_n$ can be expanded in the
basis $\{H_k\}_{k=0}^{n}$ as

\begin{eqnarray}\label{HermiteEx}
P_n(x)=H_n(x)+\alpha_{n-1}H_{n-1}(x)+\alpha_{n-2}H_{n-2}(x), \quad
\alpha_{n-k}=\frac{\int P_n(x)H_{n-k}(x)e^{-x^2}dx}{\sqrt{\int H^2_{n-k}(x) e^{-x^2}dx}}; k=1,2;
\end{eqnarray}

\noindent from where we deduce that  the  monic orthogonal polynomial $Q_{n}$ with respect to
$(\LL,\mu )$ for the  index $n$, for $n>3$, can be described as

$$Q_{n}(x)= \begin{cases}
\dsty H_n(x)+\frac{\lambda_{n}\alpha_{n-1}}{\lambda_{n-1}} H_{n-1}(x)+\frac{\lambda_{n}\alpha_{n-2}}{\lambda_{n-2}} H_{n-2}(x)+cx, & c\in \CC, \\
\dsty x.
\end{cases}$$

By a similar argument, by  expanding $P_0 ,P_1 ,P_2 ,P_3 $ in terms of $H_0 ,H_1,H_2 ,H_3 $ we have that
the solutions to \eqref{SolSet1} and \eqref{SolSet2} give that  if $n\leq 3$ then  $Q_0(x)=1, Q_1(x)=Q_2(x)=Q_3(x)=x$; therefore,  we have normality only for $n\leq 3$.

\end{example}

\begin{note}
We correct here \cite[Ex. 1]{ApLoMa02}. There, it is stated that any exactly solvable operator for which $\rho_{M}\equiv 1$
and $\rho_{k} \not\equiv 0, 0\leq k <M$ satisfies the conditions of \cite[Th. 3]{ApLoMa02}.

For the operator of Example  \ref{counterex} of the present paper, $\LL[f]=f^{\prime \prime} -2xf^{\prime}+2f$. The
function  $\rho_{0,1}$ defined in \cite[Th. 3]{ApLoMa02} (according with the notation employed in that paper) simplifies to
$\rho_0+\rho_1=2-2\,x$ and this function has a zero on $supp(\mu)$. This particular operator does not satisfy the
conditions given in \cite[Th. 3]{ApLoMa02}; therefore,  that theorem cannot be applied for exactly solvable operators
in general.
\end{note}

In order to provide a necessary and sufficient condition for the normality of an index, we introduce some auxiliary
notation and  prove some preliminaries lemmas. In the sequel, we denote by  $\Delta_{n}$  the determinant of the Hankel matrix
defined by the moments $\mu_0,\ldots, \mu_{2n}$ of the measure $\mu$. Define  $\Delta_{0,0}=\mu_0$ and denote by
$\Delta_{n,i},  0\leq i\leq n,$ the determinant of the  following matrix with  column $i+1$ deleted

$$
\left ( \begin{array}{ccc}
                                                           \mu_0      & \cdots  & \mu_{n} \\
                                                                      & \vdots  &          \\
                                                            \mu_{n-1} &  \cdots & \mu_{2n-1}
                                                       \end{array} \right ).
$$

 Consider the infinite upper triangular matrix $A=(a_{i,j})$ with entries

\begin{equation}\label{MatAdef}
a_{i,j}=\dsty \sum_{k=j-i}^{\min(M,j-1)}\rho_{k,i+k-j}\frac{(j-1)!}{(j-1-k)!},\quad i\leq j,
\end{equation}

\noindent and set $A_n=(a_{i,j})_{\begin{subarray}{c}
                           1\leq i\leq n \\
                           1\leq j \leq n
                          \end{subarray}}$.

Let $\dsty q_n(x)=\sum_{j=0}^{n}\alpha_{n,j}x^j$ be a generic polynomial of degree $n$. By $\mathfrak{a}_{n+1}=(\alpha_{n,0},\dots ,\alpha_{n,n})^t$ we denote the column vector of the coefficients of $q_n$.

\begin{lemma}\label{Iso}
Let $\mu$  be a positive Borel measure on the real line, $\left\{P_n\right\}_{n=0}^{\infty}$ the associated
sequence of monic orthogonal polynomials, and  $\dsty \LL^{(M)}=\sum_{k=0}^{M}\rho_{k}(x)\frac{d^{k}}{dx^{k}}$  an
exactly solvable operator, where $\dsty \rho_{k}(x)=\sum_{j=0}^{k}\rho_{k,j}x^j$. Then, \eqref{SolSet1}
can be expressed  as

\begin{eqnarray}
\label{syst1} A_{n+1}\,\mathfrak{a}_{n+1} &=&\lambda_n \mathfrak{b}_{n+1},
\end{eqnarray}
and  \eqref{SolSet2} can be expressed as

\begin{eqnarray}
\label{syst2} A_{n+1}\,\mathfrak{a}_{n+1}=0,
\end{eqnarray}

\noindent where  $\mathfrak{b}_{n+1}=(\beta_{n,0},\dots,\beta_{n,n})^t, \beta_{n,i}=\Delta_{n,i} \Delta_{n,n}^{-1}, 0\leq i\leq n,$
is the column vector of the coefficients of $P_n$.

\end{lemma}

\proof

Let $\dsty q_n(x)=\sum_{j=0}^{n}\alpha_{n,j}x^{j}$. We have that

\begin{eqnarray}\label{SystemEquiv}
\LL^{M}\left[\sum_{j=0}^{n}\alpha_jx^j\right]&=&\sum_{k=0}^{M}\left(\sum_{u=0}^{k}\rho_{k,u}x^u\sum_{j=0}^{n}\alpha_{n,j}x^{j-k}\frac{j!}{(j-k)!}\right)\\
\nonumber
&=&\sum_{k=0}^{M}\sum_{j=k}^{\min(k,n)}\sum_{u=0}^{k}\rho_{k,u}\frac{j!}{(j-k)!}x^{u+j-k}\alpha_{n,j}\\ \nonumber
&=&\sum_{k=0}^{M}\sum_{i=0}^{n}\sum_{j=\max(k,i)}^{\min(n,i+k)}\left(\rho_{k,i+k-j}\frac{j!}{(j-k)!}\alpha_{n,j}\right)x^{i}\\
\nonumber
&=&\sum_{i=0}^{n}\sum_{j=i}^{n}\sum_{k=j-i}^{\min(M,j)}\left(\rho_{k,i+k-j}\frac{j!}{(j-k)!}\alpha_{n,j}\right)x^{i}.
\end{eqnarray}

From Heine's formula for the monic orthogonal polynomials,  we have

\begin{eqnarray}\label{Detf}
P_n(x)= \Delta_{n,n}^{-1}\left | \begin{array}{ccc}
                                                                \mu_0             & \cdots & \mu_{n} \\
                                                                           & \vdots &             \\
                                                                \mu_{n-1}  & \cdots & \mu_{2n-1}  \\
                                                            1 &  \cdots & x^{n}
                                                       \end{array} \right |.
\end{eqnarray}

Therefore, \eqref{SolSet1} or  \eqref{SolSet2} can be expressed in matrix form as

\begin{eqnarray*}
A_{n+1}\,\mathfrak{a_{n+1}} &=&\lambda_n \mathfrak{b_{n+1}},
\end{eqnarray*}
 and

\begin{eqnarray*}
 A_{n+1}\,\mathfrak{a_{n+1}}=0,
\end{eqnarray*}
respectively.  \qed

\begin{lemma}\label{CoeffPoweru}
Let   $\dsty \LL^{(M)}=\sum_{k=0}^{M}\rho_{k}(x)\frac{d^{k}}{dx^{k}}$ be an exactly solvable operator, $\dsty
\rho_{k}(x)=\sum_{i=0}^{k}\rho_{k,i}x^i$, and $A_{n+1}$ the matrix whose entries are defined by \eqref{MatAdef}. Then,  $a_{j+1,j+1}$ is the
coefficient of $x^{j}$ of the polynomial $\dsty \LL^{(M)}[x^{j}]$.
\end{lemma}

\proof
We have
\begin{eqnarray*}
\LL^{(M)}[x^j]&=&\sum_{k=0}^{M}\rho_{k}(x)\frac{d^{k}}{dx^{k}}x^j\\
              &=&\sum_{k=0}^{M}\rho_{k}(x)\frac{j!}{(j-k)!}x^{j-k}\\
              &=&\sum_{k=0}^{\min(M,j)}\sum_{i=0}^{k}\rho_{k,i}\frac{j!}{(j-k)!}x^{i+j-k}.
\end{eqnarray*}

From this expression we obtain that if  $i=k$ then the coefficient of $x^j$ in $\LL^{(M)}[x^j]$ is

$$
\sum_{k=0}^{\min(M,j)}\rho_{k,k}\frac{j!}{(j-k)!},
$$

\noindent which corresponds to the coefficient $a_{j+1,j+1}$ in \eqref{MatAdef} of the matrix $A_{n+1}$. \qed

From the preceding lemmas we   deduce a necessary and sufficient condition for the normality of an index $n$.

\begin{theorem}\label{necsufQng}
Let $\mu$  be a positive Borel measure on the real line, $\left\{P_n\right\}_{n=0}^{\infty}$ the associated
sequence of monic orthogonal polynomials, and  $\dsty \LL^{(M)}=\sum_{k=0}^{M}\rho_{k}(x)\frac{d^{k}}{dx^{k}}$  an
exactly solvable differential operator. Then, an index $n\in \ZZ_+$ is normal if and only if either

\begin{itemize}
\item [i)] $ deg[\LL^{(M)}[x^k]]=k, \quad \forall  k:\quad 0 \leq k\leq n$,
or

\item [ii)] There exist indexes $n_1,\ldots, n_k;  0\leq n_1\leq \ldots \leq n_k\leq n,$ such that $ deg[\LL^{(M)}[x^{n_j}]]<n_j,  1\leq
    j\leq k$,

        \begin{itemize}
 \item [1)] if $k\geq 1$,  then $\{\LL^{(M)}[1],\LL^{(M)}[x],\linebreak[1]\dots ,\LL^{(M)}[x^{n_k}]\}$ has $n_k$  linearly  independent vectors,

  \item [2)] if  $n_k<n$, then  the moments of the measure $\mu$ satisfy the relation

\begin{eqnarray}
\label{NormalMom1}\sum_{j=0}^{n-n_k-1}\gamma_{n-n_k-j} \Delta_{n,n-j}\neq -\Delta_{n,n_k},
\end{eqnarray}
 where $\{\gamma_i\}_{i=1}^{n-n_k}$  are such that

\begin{eqnarray}\label{gammacoeff}
( \gamma_1,    \ldots    ,\gamma_{n-n_k})B =(-a_{n_{k}+1,n_{k}+2},  \ldots              ,-a_{n_{k}+1,n+1}),
\end{eqnarray}
and $B$ is the matrix

\begin{eqnarray*}
B=\left(
      \begin{array}{ccccc}
                                 a_{n_{k}+2,n_{k}+2}   & a_{n_{k}+2,n_{k}+3}&                    &           & \vdots      \\
                                                       &     \vdots         &                    &           &             \\
                                                       &                    &       \cdots       &   a_{n,n} &  a_{n,n+1}   \\
                                         0             &                    &      \cdots        &    0      &  a_{n+1,n+1} \\
       \end{array}
    \right).
\end{eqnarray*}

\end{itemize}

\end{itemize}
\end{theorem}

\proof We assume that $n\geq 1$, otherwise we have   that $n=0$ is normal and   there is nothing to prove. Suppose that
the index  $n$ is normal. This is equivalent to saying that the following alternatives  for \eqref{SolSet1} and
\eqref{SolSet2} hold.

\begin{itemize}
\item [a)] Equation \eqref{SolSet1} has a unique monic polynomial solution.
\item [b)] Equation \eqref{SolSet2} has a unique non zero monic polynomial solution and  $\lambda_n\neq 0$.

\item [c)] Equation  \eqref{SolSet2} has a unique non zero monic polynomial solution and $\lambda_n =0$.
\end{itemize}

If we have  alternative $a)$, then Lemma \ref{Iso} gives that this statement is equivalent to saying that  $Ker[A_{n+1}]=\{0\}$,
hence the elements of the diagonal of the matrix $A_{n+1}$ are non null. By Lemma \ref{CoeffPoweru} we obtain  that $
deg[\LL^{(M)}[x^k]]=k, \forall  k: 0 \leq k\leq n$, that is, we have $i)$.

Suppose  we have alternatives $b)$ or $c)$. We have   that there exist indexes $n_1,\ldots, n_k;  0\leq n_1\leq \ldots \leq n_k\leq n,$ such that $ deg[\LL^{(M)}[x^{n_j}]]<n_j,  1\leq
    j\leq k$.  For $k>1$, we partition  the matrix $A_{n+1} $ in blocks as
$$A_{n+1}=\left( \widetilde{B}_2,\widetilde{B}_1
    \right),
$$
where $\widetilde{B}_2 $ is the block of $A_{n+1}$ formed by its first $n_k+1$ columns,   and $\widetilde{B}_1$ has the columns $n_k+2,\ldots ,n+1$ of $A_{n+1}$. When $n_k=n$,
$$A_{n+1}=\left( \widetilde{B}_2 \right).
$$

%%%%%%%%%%%%%%%%%%%%%%%%%%%%%%%%%%%%%%%%%%%%%%%%%%%%%

Suppose  we have alternative $b)$.   If  $k\geq1$ and  $deg[\LL^{(M)}[x^{n_j}]]<n_j, \forall j=1,\ldots ,k$, by Lemma \ref{CoeffPoweru},  $a_{n_j+1,n_j+1}=0$; hence, 
$rank[\widetilde{B}_1]=n-n_k$, notice that since $\lambda_n\neq 0$ then we have that $n_k<n$, therefore the block $\widetilde{B}_1$ is not empty.  As \eqref{SolSet2} has a unique non null monic polynomial
solution and Lemma \ref{Iso} we deduce  that   $dim[Ker[A_{n+1}]]=1$; therefore,  $rank[A_{n+1}]=n$. By  denoting 
by $v_i$ the $i$--th column of $A_{n+1}$ and  by noting that   $span \left[ \{v_i\}_{i=1}^{n_k+1} \right] \bigcap\linebreak[1] span \left[ \linebreak[1]\{v_i\}_{i=n_k+2}^{n+1}\right]$ reduces to the null element, one has
\begin{eqnarray*}
n&=&rank\left[\{v_i\}_{i=1}^{n_k+1}\bigcup\{v_i\}_{i=n_k+2 }^{n+1}\right]\\ &=&
rank\left[\{v_i\}_{i=1}^{n_k+1}\right]+ rank\left[\{v_i\}_{i=n_k+2}^{n+1}\right]= rank\left[\{v_i\}_{i=1}^{n_k+1}\right]+n-n_k,
\end{eqnarray*}

%
%\begin{eqnarray*}
%n&=&rank\left[\{v_i\}_{i=1}^{n_1}\bigcup \{v_i\}_{i=n_1+1}^{n_k{\color{red}+1}}\bigcup\{v_i\}_{i=n_k+{\color{red}2}}^{n+1}\right]\\ &=&
%rank\left[\{v_i\}_{i=1}^{n_1}\right]+rank[\{v_i\}_{i=n_1+1}^{n_k+1}]+rank[\{v_i\}_{i=n_k+2}^{n+1}]=n_1+rank[\{v_i\}_{i=n_1+1}^{n_k+1}]+n-n_k,
%\end{eqnarray*}

\noindent  which implies that  $rank\left[\{v_i\}_{i=1}^{n_k+1}\right]=n_k$; which is equivalent to saying that
$\{\LL^{(M)}[1],\LL^{(M)}[x],\linebreak[1]\dots ,\LL^{(M)}[x^{n_k}]\}$ has $n_k$ linearly  independent vectors,  and we have $1)$ of $ii)$.

Consider now the statement of $b)$ that $\lambda_n\neq 0$, or equivalently, that $deg[\LL^{(M)}[x^{n}]]=n$. From Lemma \ref{CoeffPoweru} we have necessarily that $n_k<n$. Since  equation \eqref{SolSet2} has a unique non zero  monic polynomial solution and the index $n$ is normal by hypothesis, then \eqref{SolSet1} has no solution. By
Lemma \ref{Iso} we deduce that the system
 \begin{eqnarray*}
  A_{n+1}\,\mathfrak{a}_{n+1}=\lambda_n\mathfrak{b}_{n+1},
\end{eqnarray*}
 is necessarily incompatible. Let us multiply  the above relation on both sides by the matrix  
$$\Gamma=\left(
      \begin{array}{ccc}
  I_{n_k\times n_k}             &\vrule &    \theta_{n_k\times (n+1-n_k)}   \\\cline{1-3}
  \theta_{(n+1-n_k)\times n_k}  & \vrule&  \begin{array}{ccccc}
                                            1                        & \gamma_1  & \gamma_2 & \cdots  &\gamma_{n-n_k}\\
                                            0                        &   1       & 0        &  \cdots &   0           \\
                                            0                        &   0       & 1        &  \cdots &   0           \\
                                                                     &           & \vdots   &         &   \vdots       \\
                                            0                        &   0       & 0        &         &    1
      \end{array} \\
      \end{array}
    \right),
$$
where $I_{n_k\times n_k}, \theta_{(n+1-n_k)\times n_k},\theta_{n_k\times (n+1-n_k)}$ are blocks corresponding to the identity and null matrices of sizes $n_k\times n_k, (n+1-n_k)\times n_k, n_k\times (n+1-n_k)$ respectively and the $\{\gamma_i\}_{i=1}^{n-n_k}$ are as in \eqref{gammacoeff}. By noting that the row $n_k+1$ of the matrix $A_{n+1}$ is a linear combination of the rows $n_k+2,\ldots, n+1$  one has
\begin{eqnarray}\label{SolSet1aux}
  \Gamma\, A_{n+1}\,\mathfrak{a}_{n+1}=\lambda_n\,\Gamma\,\mathfrak{b}_{n+1}.
\end{eqnarray}
Hence, the system \eqref{SolSet1aux} is incompatible if and only if the component $n_k+1$ in vector
$\mathfrak{b}_{n+1}$ satisfies that
$$\sum_{j=0}^{n-n_k-1}\gamma_{n-n_k-j} \Delta_{n,n-j} \neq  -\Delta_{n,n_k},$$
and we obtain $2)$ of $ii)$.  Conversely, if $2)$ and $1)$ of $ii)$ hold then we have again the statement  $b)$ which is
equivalent to the normality of $n$.

Finally, suppose that   alternative c) is the case. Lemma \ref{Iso} gives that systems \eqref{syst1} and \eqref{syst2}
are the same. If $k=1$ then we have that $dim[Ker[A_{n+1}]]=1$ and we are done. Assume that $k>1$, then $n_k=n$ and
 $$A_{n+1}=\left(  \widetilde{B}_2 \right).
$$

 As the solution to \eqref{SolSet2} is non zero
and  unique we have that   $dim[Ker[A_{n+1}]]=1$, therefore,  $rank[A_{n+1}]=n$ which implies that
$rank[\widetilde{B}_2]=n_k$; that is, the number of linearly  independent rows of the block $\widetilde{B}_2 $ is $n_k$ and in virtue of Lemmas
\ref{Iso} and  \ref{CoeffPoweru}, this is equivalent to saying that  $\{\LL^{(M)}[1],\LL^{(M)}[x],\linebreak[1]\dots ,\LL^{(M)}[x^{n_k}]\}$  has
$n_k$ linearly independent vectors and we have $1)$ of $ii)$.  The converse implication is straightforward. \qed

It is not difficult to see that the condition $i)$ obtained in Theorem \ref{necsufQng} is equivalent to affirming  that
\begin{equation}\label{M}
\{\LL^{(M)}[1],\dots ,\LL^{(M)}[x^n]\},
\end{equation}
is linearly independent which is also equivalent to saying that this set is a Markov system.  In \cite[Th 1]{ApLoMa02} it was proved that if  \eqref{M} forms a Markov system then we have normality of an index for linear homogeneous  differential operators in general.

It seems natural  to conjecture that for a general homogeneous linear differential operator a necessary and
sufficient condition could be that, either \eqref{M} is a Markov  system or   if $k\geq 1$ then
$$\{\LL^{(M)}[1],\LL^{(M)}[x],\linebreak[1]\dots ,\LL^{(M)}[x^{n_k}]\}$$ has $n_k$  linearly  independent functions on the support of the measure $\mu$, plus some additional conditions on the moments of the measure.

\section{Existence and uniqueness of polynomial solutions of degree $n$}\label{Infinite}

An interesting phenomenon  that occurs in this type of orthogonality is the existence of operators and
measures for which the associated sequence of orthogonal polynomials reduces to a finite set. A straightforward example can
be constructed to illustrate this.

\begin{example}\label{Fistorder}[First order differential operator]
Let  $\LL [f](x)=xf^{\prime}, f \in \PP,$ and consider any positive Borel measure $\mu$ supported on
a compact subset of $\RR_+$. According to \eqref{GO}, the orthogonal polynomial $Q_n$ with respect to $(\LL,\mu) $ associated to the index
$n$ is defined by

$$\int xQ^{\prime}_{n}(x)x^{k}d\mu(x)=0, \quad  \forall k\leq n-1.$$

\noindent But this is possible if and only if  $Q^{\prime}_{n}\equiv0$. Hence the sequence
$\{Q_n\}_{n=0}^{\infty}$ reduces to a constant.

\end{example}

The preceding example  shows that the sequence of polynomials $\{Q_n\}_{n=0}^{\infty}$ orthogonal with respect to the
operator $\LL[f](x)=xf^{\prime}, f \in \PP,$ and any positive Borel measure supported on a compact subset of $\RR_+$ reduces to a constant. In this
section we analyze necessary and sufficient conditions on the pair $(\LL^{(M)},\mu)$ for the existence and
uniqueness of infinite sequences of orthogonal polynomials. We shall need the following preliminary lemma.

\begin{lemma}\label{GeneralL}
Let $\{\widehat{P}_n\}_{n=0}^{\infty}$, $deg[\widehat{P}_n]=n$, be a sequence of monic polynomials and  $\dsty
\LL^{(M)}=\sum_{k=0}^{M}\rho_{k}(x)\frac{d^{k}}{dx^{k}}$ an  exactly solvable differential operator on $\PP$. Then  the following
conditions are equivalent:

\begin{itemize}
\item [i)]  $deg[\LL^{(M)}[x^n]]=n, \dsty  \quad \forall n\geq 0$.

\item [ii)]  For every $n\in\ZZ_+$ there exists a unique monic polynomial $\widehat{Q}_n$ such that

$$\LL^{(M)}[\widehat{Q}_n]=\lambda_n\widehat{P}_n.$$

\item [iii)] $Ker[\LL^{(M)}]=\{0\}$.

\item [iv)]  $\dsty \lambda_n=\sum_{j=0}^{M} \rho_{j,j}\,
    \frac{n!}{(n-j)!}\neq 0, \quad \forall n\geq 0.$

\end{itemize}

\end{lemma}

\proof

$i)\Leftrightarrow ii)$

Suppose that  $deg[\LL^{(M)}[x^n]]=n, \forall   k\geq 0$. Then  we have that for every fixed $n_0 \geq 0$,
$\left\{\LL^{(M)}[x^n]\right\}_{n=0}^{n_0}$ is a basis of $\PP_{n_0}$. Hence, it is possible to find
$\left\{\alpha_n\right\}_{n=0}^{n_0}$ such that $\dsty \widehat{P}_{n_0}(x)=\sum_{k=0}^{n_0}\alpha_{n_0,k}\LL^{(M)}[x^k]$ and
thus,  by construction we have that there exists a unique monic polynomial $\dsty
\widehat{Q}_{n_0}(x)=\sum_{k=0}^{n_0}\alpha_kx^k$ such that $\LL^{(M)}[\widehat{Q}_n]=\lambda_n\widehat{P}_n$
holds and we get $ii)$.

Suppose now  that for some index $n_0$ we have that $deg[\LL^{(M)}[x^{n_0}]]<n_0$.  From this fact and the hypothesis that $\LL^{(M)}$
is exactly solvable, every polynomial $Q_{n_0}$ of degree less or equal to $n_0$ satisfies that
$\LL^{(M)}[\widehat{Q}_{n_{0}}]$ is a polynomial of degree  less than $n_0$ and hence it cannot satisfy   $\LL^{(M)}[\widehat{Q}_n]=\lambda_n\widehat{P}_n$. That is $ ii)\Rightarrow  i)$; therefore,
$i)\Leftrightarrow ii)$.

$ii)\Leftrightarrow iii)$.

Assume that $ii)$ holds. As $\widehat{Q}_n$ is unique, for every non negative integer we have that
$Ker[\LL^{(M)}]=\{0\}$; that is, $ii)\Rightarrow iii)$.  The converse implication is straightforward.

$i)\Leftrightarrow iv)$.

This follows from the fact that the coefficient associated to the factor $x^n$ in $\LL^{(M)}[x^n]$ is equal to
$$\dsty
\sum_{j=0}^{M} \rho_{j,j}\, \frac{k!}{(k-j)!}.$$ \qed

We characterize now the exactly solvable operators for which we can guarantee the existence and uniqueness of an infinite sequence of
orthogonal  polynomials $\{Q_n\}_{n=0}^{\infty}$, such that each  polynomial $Q_n$ has degree equal to $n$.

\begin{theorem}\label{AllSeqUniq}
Let $\mu$  be a positive Borel measure on the real line,  $\left\{P_n\right\}_{n=0}^{\infty}$ the associated
sequence of monic orthogonal polynomials, and  $\dsty \LL^{(M)}=\sum_{k=0}^{M}\rho_{k}(x)\frac{d^{k}}{dx^{k}}$  an
exactly solvable operator with $\dsty \rho_{k}(x)=\sum_{j=0}^{k}\rho_{k,j}x^j$.   Then, there exists a unique sequence
of monic polynomials $\left\{Q_{n}\right\}_{n=0}^{\infty}$, each polynomial $Q_n$ of degree equal to $n$ and
orthogonal with respect to $(\LL^{(M)},\mu)$, if and only if any of the statements of Lemma \ref{GeneralL} hold.

\end{theorem}

\proof  It follows from  Lemma \ref{GeneralL} by taking $\dsty
\left\{\widehat{P}_n \right\}_{n=0}^{\infty}=\left\{P_n \right\}_{n=0}^{\infty}$ in $ii)$. \qed

A natural question then arises. What happens if any of the conditions of Lemma \ref{GeneralL} does  not hold? It is
not difficult to see that from  the expression of $\lambda_n$ as a polynomial in $n$ given in $iv)$ of Lemma \ref{GeneralL}, only for a finite number of values this relation
will not be valid. Let us denote by $S$ the set of such indexes. In this case, it is also possible to  give necessary
and sufficient conditions on the measure $\mu$  to have an  infinite sequence $\left\{Q_{n}\right\}_{n\notin
S}$ such that $deg[Q_n]=n$. The following theorem characterizes such measures in terms of a
finite set of difference equations with given initial conditions.

\begin{theorem}\label{sufQnNS}
Let $\mu$  be a positive Borel measure on the real line,  $\left\{P_n\right\}_{n=0}^{\infty}$ the associated
sequence of monic orthogonal polynomials,   $\dsty \LL^{(M)}=\sum_{k=0}^{M}\rho_{k}(x)\frac{d^{k}}{dx^{k}}$  an
exactly solvable differential operator, and $\dsty \rho_{k}(x)=\sum_{j=0}^{k}\rho_{k,j}x^j$. Suppose  that condition
$iv)$ of Lemma \ref{GeneralL} is not satisfied  and denote by $S=\{n_1,\dots,n_k \}$ the set of   indexes for which
that condition does not hold. Then, there exists a sequence of monic polynomials $\left\{Q_{n}\right\}_{n\notin S}$
orthogonal with respect to $(\LL^{(M)},\mu)$ such that $deg[Q_n]=n$ if and only if the moments of the measure $\mu$ satisfy the system
 \begin{eqnarray}
\nonumber \sum_{v=
-M}^{n_1}\left(\sum_{k=\max(-v,0)}^{M}\sum_{i=\max(0,v)}^{\min(n_1,v+k)}(-1)^{i+n_1}\frac{n!}{(n-k)!}\Delta_{n_1,i}\rho_{k,v-i+k}\right)\mu_{n+v}&=&0,\\
\label{systQ}   & \vdots    &  \\ \nonumber
\sum_{v=-M}^{n_k}\left(\sum_{k=\max(-v,0)}^{M}\sum_{i=\max(0,v)}^{\min(n_k,v+k)}(-1)^{i+n_k}\frac{n!}{(n-k)!}\Delta_{n_k,i}\rho_{k,v-i+k}\right)\mu_{n+v}&=&0,
\end{eqnarray}

\noindent where $n_k\in S$ and $n\notin S$. Moreover, if $\mu_0,\dots ,\mu_{2n_k-1}$ are the moments of some positive measure supported on a subset of  $\RR$ satisfying  \eqref{systQ},  then for $n>n_k$  the system
\eqref{systQ} defines a linear   system  of difference equations with varying  coefficient  and  initial conditions
$\mu_0,\dots ,\mu_{2n_k-1}$.

\end{theorem}

\proof

By definition, the set $\left\{Q_{k}\right\}_{k=0}^{n}, n\notin S$  exists  if and only if for every $n\notin S$ it
is possible to find coefficients $\left\{\alpha_k\right\}_{k=0}^{n}$ such that
\begin{eqnarray*}
\dsty P_{n}(x)=\sum_{k\notin S}\alpha_k\LL^{(M)}[x^{k}].
\end{eqnarray*}

As  $\LL^{(M)}$ is exactly solvable the preceding condition is equivalent  to
\begin{eqnarray}\label{spancond1}
span[\{P_{k}\}_{k=0}^{n}]=span[\{\LL^{(M)}[x^k]\}_{k=0}^{n}], \quad k\notin S,
\end{eqnarray}

\noindent and  \eqref{spancond1} is equivalent to saying that there exist coefficients
$\left\{\beta_k\right\}_{k=0}^{n}$ such that
$$\dsty \LL^{(M)}[x^{n}]=\sum_{k\notin S}\beta_kP_{k}(x), \quad n\notin S,$$

\noindent and this condition is satisfied if and only if  $\mu$ satisfies the finite system of equations

\begin{eqnarray}
\label{sufQnNSa} \int \LL^{(M)}[x^n]P_{n_{1}}(x)d\mu(x)&=&0,\\ \nonumber \vdots  & &\\ \nonumber   \int
\LL^{(M)}[x^n]P_{n_{k}}(x)d\mu(x)&=&0,
\end{eqnarray}

\noindent for all $n\notin S$ and  $n_j\in S,  j=1,\ldots ,k$. By substituting in \eqref{sufQnNSa}  Heine's formula \eqref{Detf}
for the monic orthogonal polynomials we obtain

\begin{eqnarray*}
\nonumber \int \sum_{k=0}^{M} \frac{n!}{(n-k)!} \left | \begin{array}{ccc}
                                                                \mu_0             & \cdots & \mu_{n_{1}} \\
                                                                           & \vdots &             \\
                                                                \mu_{n_{1}-1}& \cdots  &         \mu_{2n_{1}-1} \\
                                                        \dsty \sum_{j=0}^{k}\rho_{k,j}x^{n+j-k} &  \cdots & \dsty
                                                        \sum_{j=0}^{k}\rho_{k,j}x^{n+j-k+n_{1}}
                                                       \end{array} \right | d\mu(x)&=& 0,\\
                & \vdots    &  \\
\nonumber  \int \sum_{k=0}^{M}\frac{n!}{(n-k)!}  \left | \begin{array}{ccc}
                                                                \mu_0             & \cdots & \mu_{n_{k}} \\
                                                                           & \vdots &             \\
                                                               \mu_{n_{k}-1}& \cdots  &         \mu_{2n_{k}-1} \\
                                                        \dsty \sum_{j=0}^{k}\rho_{k,j}x^{n+j-k} &  \cdots &  \dsty
                                                        \sum_{j=0}^{k}\rho_{k,j}x^{n+j-k+n_{k}}
                                                       \end{array} \right |d\mu(x) &=&0.
\end{eqnarray*}

By commuting the integral and the summation symbols, expanding the determinant  by minors and doing some change of indexes we
have 
\begin{eqnarray*}
\sum_{k=0}^{M}\sum_{i=0}^{n_1}\sum_{u=-k}^{0}(-1)^{i+n_1}\frac{n!}{(n-k)!}\Delta_{n_1,i}\rho_{k,u+k}\mu_{n+u+i}&=&0,\\
\vdots & &\\
\sum_{k=0}^{M}\sum_{i=0}^{n_k}\sum_{u=-k}^{0}(-1)^{i+n_k}\frac{n!}{(n-k)!}\Delta_{n_k,i}\rho_{k,u+k}\mu_{n+u+i}&=&0,
\end{eqnarray*}

\noindent which is equivalent to  \eqref{systQ}.

Consider now that  $\mu_0,\dots ,\mu_{2n_k-1}$ are the moments of some positive measure supported on a subset of
$\RR$ satisfying  \eqref{systQ}. It is not difficult to see that for $n>n_k$ , the  system \eqref{systQ} defines a linear system  of difference equations with  varying coefficients and with initial conditions $\mu_0,\dots ,\mu_{2n_k-1}$. \qed

For a given operator $\LL^{(M)}$, we denote the class of positive Borel measures with support contained in   $\RR$
which satisfy system \eqref{systQ} as $\Xi_{\LL^{(M)}}$. Note that in general, $\Xi_{\LL^{(M)}}$ does not necessarily reduce to the empty set, as show the following  examples.
 \begin{example}\label{FirstOp}
Consider the first order linear differential operator $\LL [f](x)=xf^{\prime}(x)-f(x),  f \in \PP$. Note that $\LL [x^n]=(n-1)x^{n}$. Hence, 
the set of indexes $n$ for which $iv)$ of Lemma \ref{GeneralL} is not fulfilled reduces to $n=1$. Then, \eqref{systQ} reads
 \begin{eqnarray*}
 \mu_0&=&c \in \RR^{+},\\
 \mu_1&=&c \in \RR,\\
 \mu_0\mu_{n+1}- \mu_1\mu_n&=&0,\quad n>1.
\end{eqnarray*}

\end{example}

\begin{example}\label{Eulerop}
Consider the Euler--Cauchy operator  $\dsty \LL^{(M)} [f](x)=\sum_{k=1}^{M}a_kx^kf^{(k)}(x)$, where  $a_k\in\RR$ are such that
the polynomial $\dsty p(n)=\sum_{k=1}^{M}\frac{n!}{(n-k)!}a_k$ does not have roots for $n>0$, notice  that $\LL^{M}
[x^n]=p(n)x^n$. If  $p(n)$ has no integer roots for $n>0$, then the  system  \eqref{systQ} reduces to,

$$\mu_n=0, \quad n\geq 1,$$

\noindent which implies that $\mu\equiv 0$. Hence, the set $\Xi_{\LL^{(M)}}$ is empty.

\end{example}

 \begin{example}\label{HermiteOp}
Let  $\dsty \LL_{H} [f](x)=f^{\prime\prime}(x)-2\,x\,f^{\prime}(x),  f \in \PP$ be the Hermite operator. Then
\eqref{systQ} is

\begin{eqnarray*}
\mu_0&=&c, \quad c\in \RR^{+},\\ \mu_1&=&0,\\ 2\mu_n-(n-1)\mu_{n-2}&=&0, \quad n\geq 2,
\end{eqnarray*}

\noindent which is the difference equation that defines the measure $c\,\mu_{H}$, where $d\mu_{H}(x)=e^{-x^2}dx$ and $c\in \RR^{+}$.

\end{example}

%$\Xi_{\LL_{S}}=\{\mu_{S}(x)\}$, $\Xi_{\LL_{L}}= \{\mu_{L}(x)\}$,

In a similar way, for the Laguerre and Jacobi operators $\LL_{L}, \LL_{(\alpha,\beta)}$, respectively, we obtain
that $\Xi_{\LL_{L}}=\{c\,\mu_{L}\}_{c\in \RR^{+}}$, $\Xi_{\LL_{(\alpha,\beta)}}=\{c\,\mu_{\alpha,\beta}\}_{c\in \RR^{+}}$, where
$d\mu_{\alpha,\beta}(x)=(1-x)^{\alpha}(1+x)^{\beta}dx, d\mu_{L}(x)=x^{\alpha}e^{-x}dx$ are the Jacobi and Laguerre
measures, respectively. As a consequence, we obtain the following corollary,

\begin{coro}
Let $\LL$ be a classical operator, i.e. Jacobi, Laguerre or Hermite and $\mu$ a positive Borel measure with
support contained in $\RR$. Then there exists an infinite sequence $\left\{Q_{n}\right\}_{n=0}^{\infty}$ of polynomials
orthogonal with respect to $(\LL,\mu)$, with $deg[Q_n]=n$  if and only if $\mu$ is one of the measures
$c\,\mu_{\alpha,\beta},c\,\mu_{L},c\,\mu_{H}; c\in \RR^+$. In such case, all the sequences  of monic orthogonal
polynomials $\left\{Q_{n}\right\}_{n=0}^{\infty}$ with respect to the pair $(\LL^{(M)},\mu)$ such that  $deg[Q_n]=n$  are of
the form $\{P_n+k_n\}_{n=0}^{\infty}$ where $\{k_n\}_{n=0}^{\infty}, k_0=0,$ is an arbitrary  sequence of
complex numbers and $\{P_n\}_{n=0}^{\infty}$ is the  sequence of monic orthogonal polynomials with respect to
$\mu$.
\end{coro}

\proof Let $\LL$ be a fixed classical operator and $\{\lambda_n\}_{n=0}^{\infty}$ the associated sequence of
eigenvalues. Then, we have  that $\lambda_n=0$ if and only if $n=0$. Hence, the system  \eqref{systQ} of Theorem \ref{sufQnNS}
reduces to a unique equation. A simple calculation yields that the moments of the measure $\mu$ coincide with the
moments of the measure of orthogonality of the sequence of eigenpolynomials of $\LL$ multiplied by a real positive constant $c$ (see Example \ref{HermiteOp} and
the comment below it). Since the moment problem for a classical measure is determinate, we obtain that $\mu$
is the measure of orthogonality of the sequence of eigenpolynomials of $\LL$ multiplied by a real positive constant $c$.

From Theorem \ref{sufQnNS} we have that for $n\geq 1$  there exists an infinite sequence
$\left\{Q_{n}\right\}_{n\geq 1}$ of polynomials orthogonal with respect to $(\LL,\mu)$ such that  $deg[Q_n]=n$ if and
only if $\mu$ is the measure of orthogonality of the sequence of eigenpolynomials of $\LL$ multiplied by a real positive constant $c$. A simple calculation shows that for $n=0$, the polynomial $Q_0=1$ satisfies the condition of orthogonality \eqref{GO} and the statement is
valid also for the sequence $\left\{Q_{n}\right\}_{n\geq 0}$.

It is not difficult to see that from  the solutions of equations \eqref{SolSet1} and \eqref{SolSet2} we obtain that all
the sequences  of monic orthogonal polynomials $\left\{Q_{n}\right\}_{n\geq 0}$ with respect to the pair
$(\LL^{(M)},\mu)$ such that  $deg[Q_n]=n$  are of the form $\{P_n+k_n\}_{n=0}^{\infty}$ where $\{k_n\}_{n=0}^{\infty}, k_0=0$ is an arbitrary  sequence of complex numbers and $\{P_n\}_{n=0}^{\infty}$ is the  sequence of monic
orthogonal polynomials with respect to $\mu$. \qed

Nevertheless, it is possible to guarantee the existence of a sequence $\left\{Q_{n}\right\}_{n> m}$, for some $m\in
\NN$ of polynomials orthogonal with respect to a classical operator for a measure $\mu$ which satisfies the
condition $d\mu^*(x)=\rho(x)d\mu(x)$ where  $\mu^*$ denotes the Jacobi, Hermite or  Laguerre measure and $\rho$ is a non
negative polynomial on the support of $\mu^*$ of degree $m$, as will be shown

\begin{lemma} \label{Theor}
Let $\LL$ be a classical operator,  $\mu$  a finite positive Borel measure on $\RR$, and $n$  a fixed  positive integer number. Then, the differential equation \eqref{SolSet1} has a unique, except an
additive constant, monic polynomial solution $Q_n$ of degree $n$ if and only if

\begin{equation}\label{classicalexistence}
\int P_n(x)d\mu^{*}(x)=0,
\end{equation}
where $P_n$ is the $n$th  monic orthogonal polynomials with respect to the measure $\mu$.
\end{lemma}

\proof

%Let $P_n$ be the $n$th monic orthogonal polynomial for $\mu$ and suppose that there exists a polynomial $Q_n$ of degree $n$, such that $\LL[Q_n]=\lambda_n\,P_n$. Then \eqref{classicalexistence} is straightforward from
%
%$$\int  \, \LL[f](x)\, d\mu^*(x)=0, \quad f \in \PP.$$

Suppose that there exists a polynomial $Q_n$ of degree $n$ such that $\LL[Q_n]=\lambda_n \,P_n$. Let us denote by $\{L_n\}$ the sequence of orthogonal polynomials with respect to the measure $\mu^*$. We have then
\begin{eqnarray}  \label{classicalexistence(9-Q)}
  Q_n(z) & = & L_n(z) + \sum_{k=0}^{n-1}a_{(n,k)} L_k(z),
  \\ P_n(z)  & = & L_n(z) + \sum_{k=0}^{n-1}b_{(n,k)}  L_k(z), \label{classicalexistence(9)}
\end{eqnarray}
where $a_{(n,k)}=\frac{\left\langle Q_n,L_k \right\rangle}{\left\langle
L_k, L_k \right\rangle}$ and $b_{(n,k)}=\frac{\left\langle L_n,
P_k \right\rangle }{\left\langle L_k, L_k \right\rangle} $.

Replacing $Q_n$ and $P_n$ in \eqref{SolSet1} by the linear combinations  (\ref{classicalexistence(9-Q)}) and  (\ref{classicalexistence(9)}),
from the linearity of $\LL[\cdot]$ and the condition that $\LL[L_n]=\lambda_n \,L_n$ we get $$b_{(n,0)}= \frac{\int L_n(x)d\mu^*(x)}{\int d\mu^{*}}= 0.$$

Conversely, assume that $P_n$ is the $n$th monic orthogonal polynomial with respect to  $\mu$ fulfilling \eqref{classicalexistence}. Let $Q_n$ be the polynomial of degree $n$ defined by
\begin{equation*}
 Q_n(z) = L_n(z) + \sum_{k=0}^{n-1}a_{(n,k)} L_k(z),
\end{equation*}
where $a_{(n,0)} = \Lambda_n$   is an arbitrary constant and $a_{(n,k)}= \frac{\lambda_{n}}{\lambda_{k}}\,\frac{\left\langle L_n,
P_k \right\rangle }{\left\langle L_k, L_k \right\rangle}$.
From the linearity of $\LL[\cdot]$ and the condition that $\LL[L_n]=\lambda_n \,L_n$ we get that $\LL[Q_n]=\lambda_n \,P_n$. \qed

As a consequence, we have

\begin{theorem} \label{coroll} Let $\LL$ be a classical operator and  $\mu$ be a finite positive Borel measure on $\RR$,   such that  $d\mu^{*}(x)=\rho(x)
d\mu(x)$, with $\rho \in {L}^2(\mu)$. Then, $m$ is the smallest natural number such that  there exists   an infinite sequence $\left\{Q_{n}\right\}_{n> m}$ of polynomials orthogonal with respect to $(\LL,\mu)$, with $deg[Q_n]=n$ if and only if $\rho$ is a polynomial of degree $m$.
\end{theorem}

\proof Suppose that $m$ is the smallest natural number such that  for each $n>m$ there exists a monic polynomial $Q_n$ of degree $n$, unique up to an additive constant and orthogonal with respect to $(\LL, \mu)$. According to  Lemma \ref{Theor}
$$\int P_n(x)d\mu^*(x)=\int P_n(x)\rho(x)d\mu(x) \left\{ \begin{array}{cc}
 =0 & \mbox{if } n>m,\\
 \neq 0 & \mbox{if } n=m.
 \end{array}
\right.$$
This is equivalent to saying that $\dsty
\rho(x)= \sum_{k=0}^{m} c_{k} P_k(x) $ with $c_m \neq 0$. The converse is straightforward.\qed

Unlike Theorem \ref{AllSeqUniq}, Theorem \ref{sufQnNS} does not guarantee the uniqueness of the sequence. A result
for the uniqueness can be obtained by fixing an adequate number of points in the complex plane, as will be shown in the next theorem. Let $\Pi=\{\pi_{m_1},\ldots,\pi_{m_n}\}$ be  a set of polynomials and $\mathcal{Z}=\{\nu_{1}, \dots ,\nu_{n}\}$ a multiset   \cite{Bli89} (that is, a set that allows repeated elements)
of points in the complex plane. We will say that $\Pi$ is an interpolating system for  $\mathcal{Z}$ if the  following relation holds
$$
 \left | \begin{array}{ccc}
\pi_{m_1}(\nu_{1})  & \cdots &\pi_{m_n}(\nu_{1})  \\ 
 & \vdots &             \\
 \pi_{m_1}(\nu_{j})  & \cdots &\pi_{m_n}(\nu_{j})  \\ 
  & \vdots &             \\
 \pi^{(m_j-1)}_{m_1}(\nu_{j})  & \cdots &\pi^{(m_j-1)}_{m_n}(\nu_{j})  \\ 
  & \vdots &             \\
 \pi_{m_n}(\nu_{n})  & \cdots  &    \pi_{m_n}(\nu_{n}) 
  \end{array} \right |\neq 0, 
$$
notice that if for some index $j$ we have $m_j$ points of the set $\mathcal{Z}$ repeated, we have  completed the $j$--th and $(j+m_j-1)$--th rows  by taking the derivatives up to order $m_j-1$. In particular, if  $\mathcal{Z}$ consists of a single point repeated $n$ times, then the above determinant coincides with the Wronskian of the system $\Pi$.
 
\begin{theorem}\label{UniSysQ}
Assume that $\mu\in \Xi_{\LL^{(M)}}$ is not empty, let $\dsty
\LL^{(M)}=\sum_{k=0}^{M}\rho_{k}(x)\frac{d^{k}}{dx^{k}}$  be an exactly solvable differential operator. Let the set
$S$ be  as defined in Theorem \ref{sufQnNS} and let us fix (allowing repeated elements) $\{\nu_{1,n}, \dots ,\nu_{j_0,n}\}$
points on the complex plane. Then,  there exists a unique monic polynomial $R_{n-n_{j_0}}$ of degree $n-n_{j_0}$ such
that
 \begin{eqnarray*}
Q_n(x)=(x-\nu_{1,n})\cdots (x-\nu_{j_0,n})R_{n-n_{j_0}}(x),
\end{eqnarray*}

\noindent is  orthogonal with respect to $(\LL^{(M)},\mu)$ provided that     $\{Q_{n_j}\}_{j=1}^{j_0}$ is an interpolating system for $\{\nu_{1,n}, \dots ,\nu_{j_0,n}\}$. Here $n\notin S,\{Q_{n_j}\}_{j=1}^{j_0}$ a basis of monic polynomial solutions to  \eqref{SolSet2}, and  $j_0$ is the largest value for which $n_{j_0}<n, n_{j_0}\in S$.
\end{theorem}

\proof

According to Theorem \ref{sufQnNS}, if $\Xi_{\LL^{(M)}}$ does not reduce to the  empty set, then there exists a sequence
$\{Q_{n}\}_{j\notin S}$ of monic polynomials orthogonal with respect to $(\LL^{(M)},\mu)$ such that  $deg[Q_{n}]=n$. Let
  $S=\{n_1,\dots, n_k\}$, notice that if $n<n_1$ then by $i)$ of Lemma \ref{necsufQng} the index $n$ is normal, hence the monic  polynomial
$Q_n$ is unique, therefore we do not have necessarily that $Q_n$ vanishes  at the points $\{\nu_{1,n}, \dots ,\nu_{j_0,n}\}$.

Assume that $n_1<n$, consider
$\{Q_{n_j}\}_{j=1}^{j_0}$ a basis of monic polynomial solutions to  \eqref{SolSet2}, and assume that
$\widehat{Q}_n$ is a monic polynomial solution of degree $n$ to \eqref{SolSet1}. Then, for a given index $n$, there
exist unique coefficients $\{\alpha_j\}_{j=1}^{j_0}$ such that any monic polynomial solution $Q_n$ of degree $n$
to  equation \eqref{SolSet1} can be expressed as
 \begin{eqnarray}\label{UniSysQ1}
Q_n(x)=\widehat{Q}_n(x)+\sum_{j=1}^{j_0}\alpha_j Q_{n_j}(x).
\end{eqnarray}

Let us consider the multiset $\{\nu_{1,n}, \dots ,\nu_{j_0,n}\}$ of $j_0$ points  on the complex
plane. To prove  the existence of a monic polynomial $R_{n-n_{j_0}}$ of degree $n-n_{j_0}$ such that
 \begin{eqnarray}\label{UniSysQ2}
Q_n(x)=(x-\nu_{1,n})\cdots (x-\nu_{j_0,n})R_{n-n_{j_0}}(x),
\end{eqnarray}
we evaluate the polynomial
$$\widehat{Q}_n(x)+\beta_1 Q_{n_1}(x)+\cdots +\beta_{j_0} Q_{n_{j_0}}(x),$$
  at the points $x=\nu_{j,n}$ and  take derivatives up to order $m_{\nu_{j,n}}-1$, where $m_{\nu_{j,n}}$ is
the number of times that the point $\nu_{j,n}$ appears  in the multiset.  We obtain that
\begin{eqnarray}
\widehat{Q}_n(\nu_1)&=&\beta_1 Q_{n_1}(\nu_1)+\cdots +\beta_{j_0} Q_{n_{j_0}}(\nu_1),\\ \nonumber & \vdots &\\
\widehat{Q}_n(\nu_{j_0})&=&\beta_1 Q_{n_1}(\nu_{j_0})+\cdots +\beta_{j_0} Q_{n_{j_0}}(\nu_{j_0}),
\end{eqnarray}
by defining  $\alpha_1,\ldots ,\alpha_{j_0} $ as the solution of the above system, we obtain the existence. The uniqueness follows immediately  by condition that the $\{Q_{n_j}\}_{j=1}^{j_0}$ is an interpolating system for $\{\nu_{1,n},
\dots ,\nu_{j_0,n}\}$. \qed

%We prove now that the polynomial $R_{n-n_{j_0}}$ is unique. Assume that there exist two different monic polynomials of
%degree $n$ that vanish at the points $\{\nu_{1,n}, \ldots ,\nu_{n_{j_0},n}\}$. Since the set
%$\{\widehat{Q}_n,Q_{n_j},\ldots, Q_{n_{j_0}}\}$       is a basis for the polynomial solutions to
%\eqref{SolSet1}, we have that
%
%\begin{eqnarray*}
%(x-\nu_{1,n})\cdots (x-\nu_{n_{j_0},n})R_{1,n-n_{j_0}}(x)&=&\widehat{Q}_n(x)+\sum_{j=1}^{n_{j_0}}\alpha_{1,j} Q_{n_j},\\
%(x-\nu_{1,n})\cdots (x-\nu_{n_{j_0},n})R_{2,n-n_{j_0}}(x)&=&\widehat{Q}_n(x)+\sum_{j=1}^{n_{j_0}}\alpha_{2,j} Q_{n_j}.
%\end{eqnarray*}
%
%
%Both expressions give
%
%$$(x-\nu_{1,n})\cdots
%(x-\nu_{n_{j_0},n})(R_{1,n-n_{j_0}}(x)-R_{2,n-n_{j_0}}(x))=\sum_{j=1}^{n_{j_0}}(\alpha_{1,j}-\alpha_{2,j}) Q_{n_j}.$$
%
%Note that the degree of the left hand side is strictly greater that  $n_{j_0}$ and the degree of the right hand side is
%at most $n_{j_0}$, which is a contradiction; that is, $R_{n-n_{j_0}}$ is unique. 

\section{Zero location of the  polynomials $Q_n$ for a subclass of exactly solvable operators}\label{Zeroloc}

In this section we study the location of the zeros of orthogonal polynomials with respect to a certain subclass of
differential operators.  We start with a discussion of the class of operators which we shall consider.

\begin{definition}\label{OpFac}
Given $M\geq 1$, we say that the linear differential operator $\LL^{(M)}$ of   $M$-th order factorizes on $\PP$
if there exist multi-indexes $(m_{1}, \ldots, m_{J})$, $(n_{1}, \ldots, n_{J})$ and polynomials
$\{\rho_{m_j}\}_{j=1}^{J}$ with $\deg[\rho_{m_{j}}]=m_{j}, j=1,\dots, J$, such that for each polynomial $\Pi_{n}\in \PP$ we have
\begin{equation}
\label{eq6} \LL^{(M)}[\Pi_{n}](z)=\left[ \rho_{m_{J}}(z)\cdots \left[ \rho_{m_{2}}(z)\left[
\rho_{m_{1}}(z)\Pi_{n}(z)\right]^{(n_{1})} \right]^{(n_{2})}\cdots \right]^{(n_{J})}.
\end{equation}
\end{definition}

If $\LL^{(M)}$ factorizes on $\PP$, we shall denote

\begin{eqnarray*}
\LL^{(n_1)}_{1}[f](z)&:=&(\rho_{m_{1}}(z)f(z))^{(n_{1})},\\ &\vdots&\\
\LL^{(n_J)}_{J}[f](z)&:=&(\rho_{m_{J}}(z)f(z))^{(n_{J})},\\
\end{eqnarray*}

\noindent and then

$$\LL^{(M)}[f]= \LL^{(n_J)}_{J}\circ\cdots\circ \LL^{(n_1)}_{1}[f].$$

We are interested in  exactly solvable operators $\LL^{(M)}$ which factorize on $\PP$, for the case in which
$\{\rho_{m_{j}}\}_{j=1}^{J}$ are polynomials with reals roots. According to Definition \ref{exactlysolv} of exactly solvable
operator, we have necessarily that

\begin{eqnarray}\label{ExactlyCond}
\dsty \sum_{k=1}^{J}m_{k}= \sum_{k=1}^{J}n_{k}=M.
\end{eqnarray}

We denote by  $C_{M}$ the convex hull of the zeros of $\dsty \prod_{i=1}^{J}\rho_{m_i}$.

If  $\LL^{(M)}$  factorizes on $\PP$, then it is not difficult to see that  $i)$ of Lemma \ref{GeneralL} is equivalent
to the condition,

\begin{eqnarray}\label{UniqueCondComp}
 \dsty \sum_{i=1}^{j}(m_{i}-n_{i})\geq 0,\quad \forall j\leq J.
\end{eqnarray}

Hence, the class of operators that factorize on $\PP$ for which there exists a unique infinite sequence
$\left\{Q_{n}\right\}_{n=0}^{\infty}$, of monic polynomials such that  $deg[Q_n]=n$, and orthogonal with respect to
$(\LL^{(M)},\mu )$, for every positive Borel measure $\mu $ supported on $\RR$, are those which satisfy
condition \eqref{UniqueCondComp}.

To locate the zeros of orthogonal polynomials with respect to operators that factorize on $\PP$ we use an integral
representation for these operators and then we apply    known theorems for  zero location  of iterated integrals of
polynomials. From the preceding discussions, it is already known that   we have cases of operators for which the
associated  sequence of orthogonal polynomials is not unique. We  will first analyze the class of operators defined by
condition \eqref{UniqueCondComp}; that is, the class for which the existence of  the full sequence
$\left\{Q_{n}\right\}_{n=0}^{\infty}$ can be guaranteed. For these operators the following integral representation
holds.

\begin{lemma}\label{rep1}
Let  $P_n$ be  the $n$-th monic orthogonal polynomial with respect to $\mu $,    $\LL^{(M)}$ is such that
factorizes on $\PP$ as $\LL^{(M)}=\LL^{(n_J)}_{J}\circ \cdots \circ \LL^{(n_1)}_{1}$ and satisfies
\eqref{UniqueCondComp}. Then, the following representation holds
 \begin{equation*}
Q_{n}= \lambda_{n} I_{1}\circ \cdots \circ I_{J}\left[P_{n} \right],
\end{equation*}
where $\dsty I_{j}$ is the integration operator, given by
 \begin{eqnarray*}
\dsty I_{j}[f](z)&=&\frac{1}{\rho_{m_{j}}(z)}\int_{z_{n_{j},j}}^{z}\int_{z_{n_{j}-1,j}}^{t_{n_{j}-1}}\cdots
\int_{z_{1,j}}^{t_{1}}f(t)dt dt_1 \cdots dt_{n_j-1},
\end{eqnarray*}
 and   $\{z_{i,j}\}_{\begin{subarray}{c}
                           i=1,\ldots,n_j\\
                           j=1,\ldots,J
                          \end{subarray}}\subset C_M$.
\end{lemma}

\proof

As \eqref{UniqueCondComp} holds, then by $ii)$ of Lemma \ref{GeneralL}  we have  $\LL^{(M)}[Q_{n}]=\lambda_{n}P_{n}$ is
solvable. Let us consider  the function
$$f(z):=\left[\rho_{m_{J}}(z) \LL^{(n_{J-1})}_{J-1}\circ \cdots \circ \LL^{(n_{1})}_{1}[Q_{n}](z)\right]^{(n_{J}-1)}.$$
Applying successively Rolle's Theorem, taking into account \eqref{UniqueCondComp}, and  that the polynomials $\rho_{m_j}$ have their zeros on
$C_M$ and are reals  we obtain  that $f$ has at least a zero $z_{1,J}$ in $C_M$. Hence, $f(z)=\lambda_n \int_{z_{1,J}}^{z}P_n(t)
dt$.

By a similar argument we will have

\begin{equation}
\label{bpq1} \rho_{m_{J}}(z)\LL^{(n_{J-1})}_{J-1}\circ \cdots \circ \LL^{(n_{1})}_1 [Q_n](z)=\lambda_n
\int_{z_{n_{J},J}}^{z}\int_{z_{n_{J}-1,J}}^{t_{n_{J}-1}}\cdots \int_{z_{1,J}}^{t_{1}} P_n(t)   dtdt_{1}\cdots
dt_{n_{J}-1},
\end{equation}
which implies that the polynomial $\int_{z_{n_{J},J}}^{z}\int_{z_{n_{J}-1,J}}^{t_{n_{J}-1}}\cdots
\int_{z_{1,J}}^{t_{1}} P_n(t)   dtdt_{1}\cdots dt_{n_{J}-1}$ is divisible by $\rho_{m_{J}}$. Therefore,  after a finite
number of steps we will have

\begin{eqnarray*}
Q_{n}(z)&=&\lambda_{n} \frac{1}{\rho_{m_{1}}(z)}\int_{z_{n_{1},1}}^{z}\int_{z_{n_{1}-1,1}}^{t_{n_{1}-1}}\cdots
\int_{z_{1,1}}^{t_{1}}   \cdots \frac{1}{\rho_{m_{J}}(z)}\\ &
&\left[\int_{z_{n_{J},J}}^{z}\int_{z_{n_{J}-1,J}}^{t_{n_{J}-1}}\cdots \int_{z_{1,J}}^{t_{1}} P_n(t)   dtdt_{1}\cdots
dt_{n_{J}-1}\right] \cdots dtdt_{1}\cdots dt_{n_{1}-1}\\ &=&\lambda_{n} I_{1}\circ \cdots  \circ I_{J}\left[P_{n}
\right](z).
\end{eqnarray*} \qed

Consider now the class of operators which do not satisfy the condition \eqref{UniqueCondComp}. A
representation similar to the one obtained in the preceding lemma can also be  given. Let us prove some preliminary lemmas.

\begin{lemma}\label{KerOpPP1}
Assume that  $\LL^{(M)}$ factorizes on $\PP$ as $\LL^{(M)}=\LL^{(n_J)}_{J}\circ \cdots \circ \LL^{(n_1)}_{1}$. Then
$Ker[\LL^{(M)}]=\{0\}$ if and only if $\LL^{(M)}[1]\neq 0$.

\end{lemma}

\proof

The implication  $Ker[\LL^{(M)}]=\{0\} \quad  \Rightarrow \quad \LL^{(M)}[1]\neq 0 $ is straightforward. Assume that
$\LL^{(M)}[1]\neq 0$. Note that
 $$deg[\LL^{(n_j)}_{j}\circ \cdots \circ \LL^{(n_1)}_{1}[1]]=\sum_{i=1}^{j}(m_{i}-n_{i}),$$
hence,  $\dsty \sum_{i=1}^{j}(m_{i}-n_{i})\geq 0, \forall j\leq J$. Therefore,  from
\eqref{UniqueCondComp} and $iii)$ of Lemma \ref{GeneralL} we obtain that $Ker[\LL^{(M)}]=\{0\}$. \qed

\begin{lemma}\label{KerOpPP2}
Let us have  that  $\LL^{(M)}$ factorizes on $\PP$ as $\LL^{(M)}=\LL^{(n_J)}_{J}\circ \cdots \circ \LL^{(n_1)}_{1}$ and denote
by $j_0$ the largest index such that $\dsty \sum_{i=1}^{j_0}(m_{i}-n_{i})<0$. Then, $Ker[\LL^{(M)}]=\{1,\dots
,x^{n^{\prime}_{j_0}}\}$, where $\dsty n^{\prime}_{j_0}=\sum_{i=1}^{j_0}(n_{i}-m_{i})-1$.

\end{lemma}

\proof
Since $\LL^{(M)}$ is a composition of operators, it is not difficult to see that if $1\leq n \leq n^{\prime}_{j_0}$ then
 $\LL^{(M)}[x^{n}]=0$. Hence, $\{1,\ldots ,x^{n^{\prime}_{j_0}}\} \subset Ker[\LL^{(M)}]$. Suppose now that $n=n^{\prime}_{j_0}+m,  m\geq 1$. Then,  we have
 $$\dsty deg[\LL^{(n_{j_0})}_{j_0}\circ \cdots \circ \LL^{(n_1)}_{1}[x^{n}]]=m-1\geq 0,$$
and thus $\LL^{(M)}[x^n]\neq 0$.  \qed

An analogue of  Lemma \ref{rep1}, for operators that do not satisfy condition \eqref{UniqueCondComp}, is

\begin{lemma}\label{rep2}
Assume that $\mu \in \Xi_{\LL^{(M)}}\neq \emptyset$,  $\LL^{(M)}$ factorizes on $\PP$ as $\LL^{(M)}=\LL^{(n_J)}_{J}\circ
\cdots \circ \LL^{(n_1)}_{1}$ and suppose that condition \eqref{UniqueCondComp} is not satisfied. Let us fix a multiset
$\mathcal{Z}=\{\nu_{1,n}, \dots ,\nu_{n^{\prime}_{j_0}+1,n}\}$ of real numbers such that $\{1,\ldots,x^{n^{\prime}_{j_0}} \}$ is an interpolating system for $\mathcal{Z}$, where   $n^{\prime}_{j_0}$ is as in Lemma \ref{KerOpPP2}, and let $Q_n$
be the monic orthogonal polynomial  with respect to $(\LL^{(M)},\mu )$ that vanishes at the points $\{\nu_{1,n}, \dots
,\nu_{n^{\prime}_{j_0}+1,n}\}$. Then, the following representation holds
\begin{equation*}
Q_{n}(x)= \lambda_{n} I_{1}\circ \cdots I_{J-1}\circ \widehat{I}_{J,n}\left[P^{(n_{j_0})}_n\right](x) , \quad
n>n^{\prime}_{j_0},
\end{equation*}
 where $\dsty I_{j}$ is the integral operator given by
 \begin{eqnarray*}
\dsty I_{j}[f](z)&=&\frac{1}{\rho_{m_{j}}(z)}\int_{z_{n_{j},j}}^{z}\int_{z_{n_{j}-1,j}}^{t_{n_{j}-1}}\cdots
\int_{z_{1,j}}^{t_{1}}f(t)dt dt_1 \cdots dt_{n_j-1},\\
\dsty \widehat{I}_{J,n}[f]&=& I_{J}\circ I_{*,n}[f], \\
I_{*,n}[f]&=&\int_{z^{*}_{n^{\prime}_{j_0}+1,J}}^{z}\int_{z^{*}_{n^{\prime}_{j_0},J}}^{t_{n^{\prime}_{j_0}}}\cdots
\int_{z^{*}_{1,J}}^{t_{1}}f(t)dt dt_1 \cdots dt_{n^{\prime}_{j_0}},
\end{eqnarray*}
 and $\{z^*_{i,J}\}_{i=1,\ldots,n^{\prime}_{j_0}+1},\{z_{i,j}\}_{\begin{subarray}{c}
                           i=1,\ldots,n_j\\
                           j=1,\ldots,J
                          \end{subarray}}\subset  C^{*}_M $, where  $C^{*}_M $ is the convex hull of the
zeros of
 $$\dsty (x-\nu_{1,n})\cdots (x-\nu_{n_{j_0},n})\left(\prod_{i=1}^{J}\rho_{m_i}(x)\right).$$

\end{lemma}

\proof

By Lemma \ref{KerOpPP2} we have that $Ker[\LL^{(M)}]=\{1,\dots ,x^{n^{\prime}_{j_0}}\}$, where $\dsty
n^{\prime}_{j_0}=\sum_{i=1}^{j_0}(n_{i}-m_{i})-1$. Set  $S=\{0,\dots ,n^{\prime}_{j_0}\}$. Theorem \ref{UniSysQ} yields that
there exists a unique monic polynomial $R_{n-n^{\prime}_{j_0}-1}$ such that if $n>n^{\prime}_{j_0}$
$$\LL^{(M)}[\Pi_{n^{\prime}_{j_0}+1}R_{n-n^{\prime}_{j_0}-1}](x)=\lambda_n P_n(x),$$
 where $\Pi_{n^{\prime}_{j_0}+1}(x)=(x-\nu_{1,n})\cdots (x-\nu_{n^{\prime}_{j_0}+1,n})$ and $P_n$ is the $n$ th  monic orthogonal
polynomial with respect to the measure $\mu $. Taking derivatives up to order $n_{j_0}$ in the above expression,
we obtain
$$\LL^{(n_J+n^{\prime}_{j_0}+1)}_{J}\circ \cdots \circ \LL^{(n_1)}_{1}[\Pi_{n^{\prime}_{j_0}+1}R_{n-n^{\prime}_{j_0}-1}](x)=\lambda_n
P^{(n^{\prime}_{j_0}+1)}_n(x)$$
or, equivalently,
 $$\widehat{\LL}[R_{n-n^{\prime}_{j_0}-1}](x)=\lambda_n P^{(n^{\prime}_{j_0}+1)}_n(x),$$
 where   $\widehat{\LL}=\widehat{\LL}^{(n_J)}_{J}\circ \LL^{(n_{J-1})}_{J-1}\circ \cdots \circ
\LL^{(n_{2})}_{2}   \circ \widehat{\LL}^{(n_1)}_{1}$,
\begin{eqnarray*}
\widehat{\LL}^{(n_J)}_{J}[f]&=&\LL^{(n_J+n^{\prime}_{j_0}+1)}_{J}[f], \quad f\in \PP, \\
\widehat{\LL}^{(n_1)}_{1}[f]&=&\LL^{(n_1)}_{1}[\Pi_{n^{\prime}_{j_0}+1} f].
\end{eqnarray*}

Since  the polynomial $R_{n-n^{\prime}_{j_0}-1}$ is unique, we deduce  that $\widehat{\LL}[1]\neq 0$, hence  Lemma \ref{KerOpPP1}
gives that  $Ker[\widehat{\LL}]= \{ 0\}$. Therefore, from the equivalence of $iii)$ of Lemma \ref{GeneralL} and the relation 
\eqref{UniqueCondComp}, $\widehat{\LL}$ satisfies $\dsty \sum_{i=1}^{j}(m_{i}-n_{i})\geq 0, \forall j\leq J$. By
Lemma \ref{rep1}, we obtain
 $$
R_{n-n^{\prime}_{j_0}-1}(x)= \lambda_{n} \widehat{I}_{1}\circ I_{2}\circ \cdots  \circ I_{J-1}\circ \widehat{I}_{J}\left[P^{(n^{\prime}_{j_0}+1)}_n
\right](x)
$$
or, equivalently,
 \begin{equation*}
Q_{n}(x)= \lambda_{n} I_{1}\circ \cdots \circ  I_{J-1} \circ \widehat{I}_{J}\left[P^{(n^{\prime}_{j_0}+1)}_n \right](x), \quad
n>n^{\prime}_{j_0},
\end{equation*}
where $\dsty I_{j}$ is the integral operator given by
 \begin{eqnarray*}
\dsty I_{j}[f](z)&=&\frac{1}{\rho_{m_{j}}(z)}\int_{z_{n_{j}},j}^{z}\int_{z_{n_{j}-1},j}^{t_{n_{j}-1}}\cdots
\int_{z_{1,j}}^{t_{1}}f(t)dt dt_1 \cdots dt_{n_j-1},\\
\dsty \widehat{I}_{1}[f]&=& \frac{1}{\Pi_{n^{\prime}_{j_0}+1}(x)}I_{1}[f], \quad f\in \PP,\\ 
\dsty \widehat{I}_{J}[f]&=& I_{J}\circ I_{*}[f], \\ 
I_{*}[f]&=&\int_{z^{*}_{n^{\prime}_{j_0}+1,J}}^{z}\int_{z^{*}_{n^{\prime}_{j_0},J}}^{t_{n^{\prime}_{j_0}}}\cdots
\int_{z^{*}_{1,J}}^{t_{1}}f(t)dt dt_1 \cdots dt_{n^{\prime}_{j_0}},
\end{eqnarray*}
  and $\{z^*_{i,J}\}_{i=1,\ldots,n^{\prime}_{j_0}+1},\{z_{i,j}\}_{\begin{subarray}{c}
                           i=1,\ldots,n_j\\
                           j=1,\ldots,J
                          \end{subarray}}\subset  C^{*}_M $, being  $C^{*}_M $ is the convex hull of the
zeros of the polynomial
$$\dsty (x-\nu_{1,n})\cdots (x-\nu_{n^{\prime}_{j_0}+1,n})\left(\prod_{i=1}^{J}\rho_{m_i}(x)\right).$$ \qed

\subsection{Zero location }

Assume that the exactly solvable operator $\LL^{(M)}$ factorizes on $\PP$ and  that there exists a unique infinite
sequence $\left\{Q_{n}\right\}_{n=0}^{\infty}$, of monic polynomials, each polynomial $Q_n$ of degree equal to $n$,
and orthogonal with respect to $(\LL^{(M)},\mu )$, for every positive Borel measure $\mu $ supported on $\RR$.
The following theorem, see \cite[Exer. 20, pag. 74]{Mar66}, and the results of the preceding section can be used now to
locate the zeros of the family $\{Q_n\}$.

\begin{theorem}\label{mardem}
If all the zeros of the $n$th degree polynomial $f$ lie in a convex region $K$ containing the point $a$, then all
the zeros of $F(z)=\int_{a}^{z}f(t)dt$ lie in the domain bounded by the envelope of all circles passing through $a$ and
having centers on the boundary of $K$.

\end{theorem}

The following lemma will be necessary for the zero location theorem.

\begin{lemma}\label{preitermult}
Let $I_j$ be the integral operator defined in Lemma \ref{rep1}. Assume that the set $\{z_{i,j}\}_{i=1}^{n_j}$ and the
zeros of the $n$th degree polynomial $\Pi_n$ lie the circle $C(0,r)$ with center in the origin and radius $r$. Then,
the zeros of $I_j[\Pi_n]$ lie in  the circle $C(0,3^{n_j}r)$.
\end{lemma}

\proof

If   $z_{i,j}$ and the zeros of the $n$th degree polynomial $\Pi_n $ lie in a circle $C(0,r)$ of radius $r$, by
Theorem \ref{mardem} the zeros of $\int_{z_{1,j}}^{z}\Pi_n(t)dt$ lie in the envelope of all the circles with center in
the boundary of $C(0,r)$ and passing through $z_{1,j}$. It is not difficult to see that this envelope and the set
$\{z_{i,j}\}_{i=2}^{n_j}$ are contained  in the circle $C(0,3r)$. Using the same argument,  we obtain that the zeros of
$I_j[\Pi_n ]$ are located in the circle $C(0,3^{n_j}r)$. \qed

Consider now the case of operators for which the full sequence of $\{Q_n\}_{n=0}^{\infty}$ exists, for every Borel
measure $\mu $ supported on a subset of $\RR$ or, equivalently,  the operators for which this can be guaranteed are
those which satisfy the condition \eqref{UniqueCondComp}. We have then,

\begin{theorem}\label{ZeroLocUniq}
Let $\LL^{(M)}$ be an exactly solvable operator that factorizes on $\PP$ satisfying the condition
\eqref{UniqueCondComp} and  $\mu $ a positive Borel measure supported in $[-1,1]$. Then, the zeros of the sequence
$\left\{Q_{n}\right\}_{n=0}^{\infty}$, of monic polynomials orthogonal with respect to $(\LL^{(M)},\mu )$  are
located in a circle of radius $R$, where $R=3^Md$, with $\dsty d=\max\{1,\sup_{z\in C_M}|z|\}$.
\end{theorem}

\proof

As $\{P_n\}_{n=0}^{\infty}$ is the sequence of orthogonal with respect to  $\mu $, their zeros are  in $[-1,1]$. Notice
that the interval $C_M$ and the zeros of the sequence $\{P_n\}_{n=0}^{\infty}$ are contained in a circle with center at the origin and
radius $\dsty d=\max\{1,\sup_{z\in C_M}|z|\}$. From Lemma \ref{rep1} we have that $Q_{n}$ can be represented as $Q_{n}(z)= \lambda_{n} I_{1}\circ \cdots \circ I_{J}\left[P_{n} \right](z)$. Applying successively  Lemma \ref{preitermult} we obtain that the zeros are located in a circle of radius $R$, where
$R=3^Md$, with $\dsty d=\max\{1,\sup_{z\in C_M}|z|\}$.\qed

Consider now the class of operators which do not satisfy the condition \eqref{UniqueCondComp}. In this case the
associated  sequence of orthogonal polynomials is not unique, nevertheless, in Theorem \ref{UniSysQ} it was shown that
if we fix an adequate number of points we can define a unique infinite sequence of orthogonal polynomials. We have

\begin{theorem}\label{ZeroLocNonUniq}
Let $\LL^{(M)}$ be an exactly solvable operator that factorizes on $\PP$ and assume that  condition
\eqref{UniqueCondComp} is not satisfied, $\mu \in \Xi_{\LL^{(M)}}\neq \emptyset $ such that $supp(\mu )\subset [-1,1]$ and consider a sequence of multisets $\{\nu_{1,n}, \ldots ,\nu_{n^{\prime}_{j_0}+1,n}\}$,
satisfying the hypothesis of  Lemma \ref{KerOpPP2}. Then, the zeros of the sequence
$\left\{Q_{n}\right\}_{n=n_0+1}^{\infty}$ of monic orthogonal  polynomials with respect to $(\LL^{(M)},\mu )$ such that $Q_{n}(\nu_{j,n})=0, 1\leq j\leq  n^{\prime}_{j_0}+1$ are
located in a circle of radius $R$, where $R=3^Md$, with $\dsty d=\max\{1,\sup_{z\in C^{*}_M}|z|\},
C^{*}_M=\sup_{n} \bigcup_{j=0}^{n} C^{*}_{M,j}$, being $C^{*}_{M,n} $ is the convex hull of the zeros of $\dsty (x-\nu_{1,n})\cdots
(x-\nu_{n^{\prime}_{j_0}+1,n})\left(\prod_{i=1}^{J}\rho_{m_i}(x)\right)$.
\end{theorem}

\proof

By hypothesis, the zeros of the sequence $\{P_n\}_{n=0}^{\infty}$ are contained in $[-1,1]$. According to Theorem \ref{UniSysQ} there
exists a unique sequence $\left\{Q_{n}\right\}_{n=n^{\prime}_0+1}^{\infty}$ of monic orthogonal  polynomials with respect to
$(\LL^{(M)},\mu )$. By Lemma \ref{rep2},
$$Q_{n}(x)= \lambda_{n} I_{1}\circ \cdots \circ \widehat{I}_{J}\left[P^{(n^{\prime}_{j_0}+1)}_n\right](x), \quad n>n^{\prime}_{j_0},$$

\noindent where $\dsty I_{j}$ is the integral operator given by
 \begin{eqnarray*}
I_{j}[f](z)&=&\frac{1}{\rho_{m_{j}}(z)}\int_{z_{n_{j},j}}^{z}\int_{z_{n_{j}-1,j}}^{t_{n_{j}-1}}\cdots
\int_{z_{1,j}}^{t_{1}}f(t)dt dt_1 \cdots dt_{n_j-1},\\
\dsty \widehat{I}_{J}[f]&=& I_{J}\circ I_{*}[f],\\
I_{*}[f]&=&\int_{z^{*}_{n^{\prime}_{j_0}+1,J}}^{z}\int_{z^{*}_{n^{\prime}_{j_0},J}}^{t_{n^{\prime}_{j_0}}}\cdots
\int_{z^{*}_{1,J}}^{t_{1}}f(t)dt dt_1 \cdots dt_{n^{\prime}_{j_0}},
\end{eqnarray*}
  and $\{z^*_{i,J}\}_{i=1,\ldots,n^{\prime}_{j_0}+1},\{z_{i,j}\}_{\begin{subarray}{c}
                           i=1,\ldots,n_j\\
                           j=1,\ldots,J
                          \end{subarray}}\subset  C^{*}_{M,n} $, being  $C^{*}_{M,n} $ is the convex hull of the
zeros of the polynomial
$$\dsty (x-\nu_{1,n})\cdots (x-\nu_{n^{\prime}_{j_0}+1,n})\left(\prod_{i=1}^{J}\rho_{m_i}(x)\right).$$

Note that the set $C^{*}_{M,n}$ and the zeros of $\{P_n\}_{n\geq0}$ are contained in a circle with center at the origin and
radius $\dsty d=\max\{1,\sup_{z\in C^{*}_{M,n}}|z|\}$. Applying successively  Lemma \ref{preitermult} we obtain that
the zeros of $Q_n$ are located in a circle of radius $R_n$, where $R_n=3^Md_n$, with $\dsty d_n=\max\{1,\sup_{z\in
C^{*}_{M,n}}|z|\}$. Hence, the zeros of the  full sequence are located in a circle of radius
$R=3^Md$, with $\dsty d=\max\{1,\sup_{z\in C^{*}_M}|z|\},  C^{*}_M=\sup_{n} \bigcup_{j=0}^{n} C^{*}_{M,j}$. \qed

\section{The polar polynomials case}\label{Polar}

In this section we study analytic properties of the polar polynomials, already introduced in Section \ref{examples}.
Let us denote by   $\dsty d\mu_{T}(x)= \frac{1}{\sqrt{1-x^2}}dx$ the first kind Chebyshev measure and by $T_n$ the
$n$-th Chebyshev monic polynomial of the first kind. We shall study these polynomials for the   class of finite positive Borel measures
on $[-1,1]$ defined as  $\dsty d\mu (x)= \frac{d\mu_{T}(x)}{\rho(x)}$ with  $\dsty \rho(z)= r \prod_{i=1}^{m} (z-\nu_i)$ a
non negative polynomial on $[-1,1]$. Denote  by $\mathcal{P}_m(\mu_{T})$ this class of measures. This
complements the study carried out in \cite{BePiMaUr11} where the measure  $\mu $ is the Gegenbauer measure. We   obtain a curve
which contains the accumulation points of the zeros of these polynomials and a formula for the strong asymptotic
behavior of these polynomials in $\CC\setminus [-1,1]$.

\subsection{Strong asymptotic behavior and zero location}

We recall that a measure supported on $[-1,1]$ is in the Szeg\"{o} class $\mathfrak{S}$ if its absolutely continuous part
$\mu^{\prime} $ satisfies

$$\int_{-1}^{1}\frac{\log \mu^{\prime}(x) }{\sqrt{1-x^2}}dx>-\infty.$$

The asymptotic properties of orthogonal polynomials with respect to a measure supported on $[-1,1]$ in the Szeg\"{o}
class can be described by means of the Szeg\"{o} function $D(\mu ,z)$, cf. \cite[\S 6.1]{Nev79}.

\begin{definition}\label{szegof}
Let $\mu\in \mathfrak{S}$, then the Szeg\"{o} function $D(d\mu,z)$ is defined by

$$D(\mu(x),z)=\exp\left[\frac{1}{4\pi}\int_{-\pi}^{\pi}\log{\mu^{\prime}(\cos(t))}\frac{1+ze^{-\imath
t}}{1-ze^{-\imath t}}dt \right]$$
for $|z|<1$.
\end{definition}

Let $\varphi(z)=z+\sqrt{z^2-1}$,  we take the branch of $\sqrt{z^2-1}$ for which $|\varphi(z)|>1$ whenever $z \in \CC
\setminus [-1,1]$.

It is well known that orthogonal polynomials with respect to a measure which belongs to the Szeg\"{o} class have the
following outer strong asymptotic behavior, cf. \cite[\S 6.1  Lemma 18, page 67]{Nev79},

\begin{lemma}\label{StrAsyPol}
 Let $\mu $ be a positive Borel  measure supported on $[-1,1]$, $P_n$ the $n$th monic orthogonal polynomials associated to $\mu $. Then

\begin{equation*}
\frac{\delta_nP_n(z)}{\varphi(z)^n}\rightrightarrows\frac{1}{\sqrt{2\pi}}\left(D(\sqrt{1-x^2}d\mu(x),\varphi(z)^{-1})\right)^{-1},
\end{equation*}
uniformly on compact subsets of $ \overline{\CC}  \setminus [-1,1]$, where $\delta_n$ denotes the leading coefficient of the corresponding orthonormal polynomial  of degree $n$.

\end{lemma}

The next lemmas are essential in the proof of the main theorem of this section

\begin{lemma}\label{coefcalc}
Suppose that $\mu  \in \mathcal{P}_m(\mu_{T} )$, then

\begin{equation} \label{lmt}
P_n(z)=\sum_{k=0}^{m}b_{n,n-k}T_{n-k}(z), \quad
b_{n,n-k}=\frac{\int_{-1}^{1}P_n(x)T_{n-k}(x)d\mu_{T}(x)}{\int_{-1}^{1}T_{n-k}^2(x)d\mu_{T}(x)},
\end{equation}

\noindent where $P_n,T_{n}$ are the monic orthogonal polynomials associated to the measures $\mu ,\mu_{T} $,
respectively, and   the $b_{n-k,k}$ satisfy
 \begin{equation} \label{bn}
\lim_{n\to\infty}b_{n,n-k}=2^{m-k}a_{m-k}, \quad 0\leq k\leq m,
\end{equation}
 where $\dsty a_{k}=(-1)^{k}\sum_{1\leq \nu_1<\dots <\nu_k\leq m} u^{-1}_{\nu_1} \dots u^{-1}_{\nu_k}, \quad
u_{\nu_k}=\varphi(\nu_k)$.

\end{lemma}

\proof If $\mu \in \mathcal{P}_m(\mu_{T} )$ then  $d\mu_{T}(x)=\rho(x)d\mu(x)$, where  $\dsty \rho(x)=
\prod_{i=1}^{m} (x-\nu_i)$ is nonnegative on $[-1,1]$. Therefore

$$P_n(z)=\sum_{k=0}^{m}b_{n,n-k}T_{n-k}(z), \quad
b_{n,n-k}=\frac{\int_{-1}^{1}P_n(x)T_{n-k}(x)d\mu_{T}(x)}{\int_{-1}^{1}T_{n-k}^2(x)d\mu_{T}(x)}.$$

\noindent Hence, if $z=\frac{1}{2}(u+u^{-1})$ then $\dsty T_{n-k}(z)=\frac{u^{n-k}+u^{k-n}}{2^{n-k}}$, and

\begin{equation} \label{Tch1}
\frac{2^nP_n(z)}{u^n}= \sum_{k=0}^{m}2^kb_{n,n-k}u^{-k}+\frac{1}{u^{2n}}\sum_{k=0}^{m}2^kb_{n,n-k}u^k.
\end{equation}

From   \cite[\S 6.1 theorem 26]{Nev79} and Definition \ref{szegof}, we have

\begin{equation}\label{MCoeff}
\lim_{n\to\infty} \delta_n 2^{-n}=\frac{1}{\sqrt{2\pi}}D(\rho(x),0).
\end{equation}

From   Lemma \ref{StrAsyPol} and \eqref{MCoeff}, we obtain

\begin{equation} \label{Tch2}
\frac{2^nP_n(z)}{u^n} \rightrightarrows  (D(\rho(x),0))^{-1}\left(D(\frac{1}{\rho(x)},\varphi(z)^{-1})\right)^{-1}=
(D(\rho(x),0))^{-1} D(\rho(x),\varphi(z)^{-1}),
\end{equation}

\noindent uniformly on closed subsets of $\overline{\CC }\setminus [-1,1]$. By \cite[\S 6.1 Lemma 19]{Nev79} and Definition
\ref{szegof}

\begin{eqnarray*}
D(\rho,\varphi(z)^{-1}) & = &  2^m \; \exp\left( \frac{1}{2\pi} \int_{-1}^{1} \frac{\log(\rho(t))}{\sqrt{1-t^2}}
dt\right) \prod_{k=1}^{m} \frac{z-\nu_k}{\varphi(z)-\varphi(\nu_k)}, \\ D(\rho,0)& =
&\exp\left(\frac{1}{2\pi}\int_{-1}^{1} \frac{\log(\rho(t))}{\sqrt{1-t^2}} dt \right) .
\end{eqnarray*}

\noindent Hence, if $z=\frac{1}{2}(u+u^{-1})$, using Vieta's formula, we have that  the following identity holds

\begin{equation} \label{Tch3}
2^m \prod_{k=1}^{m} \frac{z-\nu_k}{\varphi(z)-\varphi(\nu_k)}=2^m
\prod_{k=1}^{m}\left(1-\frac{1}{uu_{\nu_k}}\right)=\sum_{k=0}^{m} 2^m a_{k}u^{k-m},
\end{equation}

\noindent where $\dsty a_{k}=(-1)^{k}\sum_{1\leq \nu_1<\dots <\nu_k\leq m} u^{-1}_{\nu_1} \dots u^{-1}_{\nu_k},
u_{\nu_k}=\varphi(\nu_k)$.

From \eqref{Tch1},\eqref{Tch2}, and \eqref{Tch3},

$$\sum_{k=0}^{m}(2^kb_{n,n-k}-2^ma_{m-k})u^{-k}+\frac{1}{u^{2n}}\sum_{k=0}^{m}2^kb_{n,n-k}u^k\rightrightarrows 0,$$

\noindent uniformly on compact subsets of $\overline{\CC} \setminus [-1,1]$. Therefore, $\dsty \lim_{n\to\infty}b_{n,n-k}=2^{m-k}a_{m-k}, \quad 0\leq k\leq m.$ \qed

\begin{lemma}\label{eqprim}
Suppose that $\mu  \in \mathcal{P}_m(\mu_{T})$. If $K$ is a compact subset of $\overline{\CC}\setminus [-1,1]$ and
$\zeta\in \CC\setminus [-1,1]$ then

$$(z-\zeta)Q_{n}(z)=\frac{u^{n+1}}{2^{n+1}(n+1)}\Psi_n(u) -\frac{u_{\zeta}^{n+1}}{2^{n+1}(n+1)}\Psi_n(u_{\zeta}),\quad
z=\frac{1}{2}(u+u^{-1}),  \quad  |u|>1 ,$$

$$\Psi_n(u)\rightrightarrows  \left(1-\frac{1}{u^{2}}\right) \left(D(\rho,0)\right)^{-1}D(\rho,u^{-1}), \quad
u_{\zeta}=\varphi(\zeta). $$

\end{lemma}

\proof

From the definition of the polynomials $Q_{n}$, we have

\begin{eqnarray}\label{defQp}
\dsty (z-\zeta)Q_{n}(z)=\int_{\zeta}^{z}P_n(t)dt.
\end{eqnarray}

\noindent From  \eqref{defQp} and \eqref{lmt}, it follows that

$$\dsty (z-\zeta)Q_{n}(z)=\int_{\zeta}^{z}P_n(t)dt=\int_{\zeta}^{z}\sum_{k=0}^{m}b_{n,n-k}T_{n-k}(t)dt.$$

\noindent Making the change of variables $\dsty t=\frac{u+u^{-1}}{2}$, we obtain

\begin{eqnarray}\label{defQp1}\dsty
(z-\zeta)Q_{n}(z)=\int_{\zeta}^{z}P_n(t)dt=\int_{\varphi(\zeta)}^{\varphi(z)}\sum_{k=0}^{m}b_{n,n-k}T_{n-k}\left(\frac{u+u^{-1}}{2}\right)\left(\frac{1}{2}-\frac{1}{2u^{2}}\right)du
\end{eqnarray}

Taking into account that for $n>m+1$

$$ \int T_{n-k}\left(\frac{u+u^{-1}}{2}\right)\left(\frac{1}{2}-\frac{1}{2u^{2}}\right)du=\int
\left(\frac{u^{n-k}+u^{-n+k}}{2^{n-k}}\right)\left(\frac{1}{2}-\frac{1}{2u^{2}}\right)du=$$
\begin{eqnarray}
 \label{TchebPrim} \dsty\frac{u^{n+1}}{2^{n+1}(n+1)}\left(g_{n-k}(u)-g_{n-2-k}(u)\right)+C,
\end{eqnarray}

\noindent where
 $$ \dsty
 g_{n-k}(u)=\frac{\dsty\frac{1}{2^{n-k}}\left(\frac{(n+1)}{(n-k+1)}u^{-k+n}+\frac{(n+1)}{(-n+k+1)}u^{k-n}\right)}{(\frac{u}{2})^n}.$$

 Hence, if we denote  $\dsty \Psi_n(u)=\sum_{k=0}^{m} b_{n,n-k}(g_{n-k}(u)-g_{n-2-k}(u))$,  from  \eqref{defQp1} and \eqref{TchebPrim}, we obtain

$$(z-\zeta)Q_{n}(z)=\left .\frac{u^{n+1}}{2^{n+1}(n+1)}\Psi_n(u)\right|_{\varphi(\zeta)}^{\varphi(z)}.$$

Using Lemma \ref{lmt},    \eqref{Tch2}  of Lemma \ref{coefcalc}, and taking into account that $\dsty \lim_{n\to \infty}
\frac{(n+1)}{(n-k-1)}=1,   0\leq k\leq m$, we obtain that

$$ \Psi_n(u)\rightrightarrows  \left(1-\frac{1}{u^{2}}\right) (D(\rho(x),0))^{-1} D(\rho(x),u^{-1}),  \quad  |u|>1.$$ \qed

\begin{theorem}\label{PiBeUrext}
Suppose that $\mu  \in \mathcal{P}_m(\mu_{T})$, where  $\mu_{T}$ is the first kind Chebyshev measure. Then the
accumulation points of zeros of $\{Q_{n}\}_{n=0}^{\infty}$ are located on the set  $E=\mathcal{E}(\zeta) \bigcup [-1,1]$, where
$\mathcal{E}(\zeta)$ is the ellipse
\begin{equation}
    \mathcal{E}(\zeta):= \left\{ z \in \CC : z= \cosh(\eta_{\zeta}+i\theta), 0 \leq \theta  < 2 \pi\right\},
\end{equation} and $\eta_{\zeta}:=\ln |\varphi(\zeta)|=\ln |\zeta+\sqrt{\zeta^2-1}|$.
    If $\dsty \delta(\zeta)>2$ then $E=\mathcal{E}(\zeta)$.

\end{theorem}

\proof From Lemma \ref{eqprim}, the zeros of $Q_n$ satisfy that

\begin{equation}\label{Tchzeroprim}
\left|\Psi_n(u)\frac{u^{n+1}}{2^{n+1}(n+1)}\right|^{\frac{1}{n}}=\left|\Psi_n(u_{\zeta})\frac{u_{\zeta}^{n+1}}{2^{n+1}(n+1)}\right|^{\frac{1}{n}},\quad
z=\frac{1}{2}(u+u^{-1}),  \quad  |u|>1,
\end{equation}

\noindent and from the definition of the function $\Psi_n$, we have that

$$
\lim_{n\to\infty}\left|\Psi_n(u)\right|^{\frac{1}{n}}=1, \quad |u|>1.
$$

\noindent Therefore, taking limits on both sides of \eqref{Tchzeroprim}, we deduce that the zeros of $Q_n$ can not
accumulate outside the set

$$\left\{ z \in \CC: |z+\sqrt{z^2-1}|=e^{\eta_{\zeta}}\right\} \, \bigcup \, [-1,1]\, \bigcup\{\zeta\}. $$

Hence, if $z$ is an accumulation point of zeros of the polynomials $Q_n$, we have that $z+\sqrt{z^2-1}= \, e^{\eta_{\zeta}+i\theta}$ and $z-\sqrt{z^2-1} = e^{-(\eta_{\zeta}+i\theta)}$ for $0 \leq
\theta  < 2 \pi$, and $2z = e^{\eta_{\zeta}+i\theta}+ e^{-(\eta_{\zeta}+i\theta)}.$ \qed
\section{Concluding remarks}\label{Cr}

Theorem \ref{necsufQng} can be extended to operators with polynomial coefficients in general, giving results   for example, for   the case of Heine--Stieltjes operators. For this, we can use the same
technique of  considering the expression of the matrix $A_{n+1}$ for these operators, which is not difficult to
construct.  It would be interesting to obtain  results similar to Theorems  \ref{sufQnNS} and \ref{UniSysQ}   to classify the measures for which it is possible to ensure the existence and uniqueness, in some sense, of orthogonal polynomials with respect to $(\LL,\mu )$ for $n>m,  m\in\ZZ_+$. It would also be  of interest   to obtain results on the zero location and asymptotic behavior for   more general classes of operators as well as to consider  more general relations of orthogonality, for instance,
\begin{eqnarray*}
\dsty \int \LL_{0}^{(M)}[Q_n(x)]P(x)d\mu(x)+ \dots +\int
\LL_{k}^{(M+k)}[Q_n(x)]P^{(k)}(x)d\mu_k(x) = 0,
\end{eqnarray*}
\noindent for any polynomial $P$ such that $deg[P]\leq n-1$, where the $\left\{\LL_{j}^{(M+k)}\right\}_{j=0}^{k}$ are linear homogeneous differential differential operators
with coefficients satisfying conditions analogous to \eqref{GO}.

\section{Acknowledgments. }
The author would like to thank  the comments and suggestions made by Professors  Andrei Mart\'inez Finkelshtein,
H\'ector Pijeira Cabrera, and an anonymous reviewer which helped improve the manuscript.

%% The Appendices part is started with the command \appendix; %% appendix sections are then done as normal sections %% \appendix

%% \section{} %% \label{}

%% References %% %% Following citation commands can be used in the body text: %% Usage of \cite is as follows: %%   \cite{key}          ==>>  [#] %%   \cite[chap. 2]{key} ==>>  [#, chap. 2] %%   \citet{key}         ==>>  Author [#]


%\bibliographystyle{model1a-num-names} %\bibliography{<your-bib-database>}

%% Authors are advised to submit their bibtex database files. They are %% requested to list a bibtex style file in the manuscript if they do %% not want to use model1a-num-names.bst.

\end{document}